\documentclass[10pt,oneside]{article}
\usepackage{graphicx} % to include eps graphics
\usepackage{epstopdf}
\usepackage{verbatim}
\usepackage{amsmath} % need for subequations
\usepackage{amsthm}
\usepackage{amssymb}
\usepackage{wasysym}
\usepackage[overload]{textcase}
\usepackage{vmargin}
\usepackage{dsfont}

\newtheorem{theorem}{Theorem}%[section]
\newtheorem{lemma}{Lemma}%[section]
\theoremstyle{definition} 
\newtheorem{example}{Example}%[section]
\newcommand{\pEnclose}[1]{\left( #1 \right)} 
\newcommand{\bEnclose}[1]{\left[ #1 \right]} 
\newcommand{\cEnclose}[1]{\{ #1 \}}

\newcommand{\dist}{d}

\newcommand{\ooint}[2]{\left( #1, #2 \right)} 
\newcommand{\ocint}[2]{\left( #1, #2 \right]} 
\newcommand{\ccint}[2]{\left[ #1, #2 \right]} 

\newcommand{\approaches}{\rightarrow}

\newcommand{\pr}[1]{\text{Prob} \left( #1 \right)} 
\newcommand{\cPr}[2]{\text{Prob} \left( #1  |  #2 \right)}
\newcommand{\E}[1]{E \left[ #1 \right]} 
\newcommand{\cE}[2]{E \left[ #1  |  #2 \right]} 
\newcommand{\transpose}{\text{\sc{T}}}
\newcommand{\id}{I}
\newcommand{\ones}{\mathds{1}}
\newcommand{\graph}{\mathcal{G}}
\newcommand{\vSet}{\mathcal{V}}
\newcommand{\eSet}{\mathcal{E}}
\newcommand{\nNodes}{N}
\newcommand{\leadsTo}{\to}
\newcommand{\outDeg}{n}
\newcommand{\leafNode}{\ell}
\newcommand{\adjM}{A}
\newcommand{\adjMEV}{\lambda}

\newcommand{\abs}[1]{\left| #1 \right|}
\newcommand{\ith}[1]{#1^{\text{th}}} 
\newcommand{\stratSet}{\Sigma}
\newcommand{\cProbSymb}{p}
\newcommand{\cChoiceSymb}{C}
\newcommand{\gProbSymb}{q}
\newcommand{\wagerSymb}{w}
\newcommand{\gGuessSymb}{G}
\newcommand{\gParam}{\beta}
\newcommand{\fort}{F}
\newcommand{\dFort}{D}
\newcommand{\valSymb}{v}
\newcommand{\hMean}{H}
\newcommand{\tmNdx}{t}
\newcommand{\tmNds}{s}
\newcommand{\dFactor}{d}
\newcommand{\pos}{X}
\newcommand{\propESymb}{m}
\newcommand{\propMEV}{r}
\newcommand{\lEV}{y}
\newcommand{\rEV}{x}
\newcommand{\invMsr}{\mu}
\newcommand{\sTime}{T}
\newcommand{\sTimeProb}{q}
\newcommand{\eSTime}{\tau}
\newcommand{\tProb}{\rho}
\newcommand{\propMntnt}{A}
\newcommand{\propMntt}{B}

\newcommand{\cStratSet}{\stratSet_c}
\newcommand{\cProb}{\cProbSymb}
\newcommand{\cProbM}{\MakeTextUppercase{\cProb}}
\newcommand{\cChoice}{\cChoiceSymb}
\newcommand{\gStratSet}{\stratSet_{g}}
\newcommand{\gProb}{\gProbSymb}
\newcommand{\gGuess}{\gGuessSymb}
\newcommand{\gWager}{\MakeTextUppercase{\wagerSymb}}
\newcommand{\wager}{\wagerSymb}
\newcommand{\wagerc}{\wager_{c}}
\newcommand{\val}{\valSymb}
\newcommand{\valM}{\MakeTextUppercase{\valSymb}}
\newcommand{\rValu}{\val^{-1}}
\newcommand{\rVal}{u}
\newcommand{\rValD}{\tilde{\rVal}}
\newcommand{\propE}{\propESymb}
\newcommand{\propM}{\MakeTextUppercase{\propE}}
\newcommand{\ntl}{\textsc{nt}}
\newcommand{\tl}{\textsc{t}}
\newcommand{\ssE}[1]{E_{\invMsr} \left[ #1 \right]}
\newcommand{\sscE}[2]{E_{\invMsr} \left[ #1 | #2 \right]}

\newcommand{\expSymb}{$\blacksquare$}

\begin{document}
\title{A Path Guessing Game with Wagering}
\author{Marcus Pendergrass \\mpendergrass@hsc.edu}
\date{}
\maketitle

\abstract{We consider a two-player game in which the first player (the Guesser) tries to guess, edge-by-edge, the path that second player (the Chooser) takes through a directed graph.  At each step, the Guesser makes a wager as to the correctness of her guess, and receives a payoff proportional to her wager if she is correct.  We derive optimal strategies for both players for various classes of graphs, and describe the Markov-chain dynamics of the game under optimal play.  These results are applied to the infinite-duration Lying Oracle Game, in which the Guesser must use information provided by an unreliable Oracle to predict the outcome of a coin toss. }



%\newpage

\section{Introduction} \label{s:introduction}
In this paper we study a two-player game in which the first player (the Guesser) tries to guess, edge-by-edge, the path that the second player (the Chooser) takes through a directed graph.  At each step, the Guesser makes a wager as to the correctness of her guess, and receives a payoff proportional to her wager if she is correct.  Optimal strategies for both players are derived for various classes of graphs, and the Markov-chain dynamics of the game are analyzed.  

The Path Guessing Game studied here is a generalization of the \emph{Lying Oracle Game} \cite{Koether, Koether2}.  
In the Lying Oracle Game, an Oracle makes a sequence of $n$ statements, at most $k$ of which can be lies, and a Guesser makes bets on whether the Oracle's next statement will be a lie or not.  
We will see that the Lying Oracle Game is equivalent to our Path Guessing Game on a certain graph whose maximum out-degree is two. 
Ravikumar in \cite{Ravikumar} demonstrates a reciprocal relationship between the Lying Oracle problem and the continuous version of \emph{Ulam's Liar Game}. In that game, a Questioner tries to find a subset of smallest measure that contains an unknown number in $\ccint{0}{1}$ by asking a Responder $n$ questions about the number's location in the interval.  Again the Responder may lie up to $k$ times.
Under optimal play, the measure of the Questioner's subset is the reciprocal of the Bettor's fortune in the Lying Oracle Game.
In \cite{Rivest} Rivest et al. use this game to analyze binary search in the presence of errors.

In addition to its intrinsic interest, the Path Guessing Game provides a context in which new questions about the Lying Oracle Game can be asked and answered.
For instance, what are the optimal strategies for the \emph{infinite-duration} Lying Oracle Game, in which no block of $n$ statements can contain more than $k$ lies?  Questions of this sort are taken up in the last section of this paper, after a general analysis of the Path Guessing Game has been carried out.

\subsection{Description of the Game} \label{ss:description}

To describe the Path Guessing Game precisely, let $\graph = \pEnclose{\vSet, \eSet}$ be a directed graph with vertex set $\vSet$ and edge set $\eSet$.  Call a vertex $j \in \vSet$ a \emph{terminal node} if it has out-degree zero, and assume that each terminal node $j$ has been assigned a positive \emph{value} $\val_{j}$.  Both players know these values.  The Guesser and the Chooser are initially located at the same non-terminal node $i$ of $\graph$, and the Guesser has an initial fortune of one dollar.  From the current node, the Chooser is going to select one of the adjacent nodes to move to next.  But before he does so, the Guesser writes down a prediction as to which node the Chooser will select.  In addition, the Guesser makes a wager with the Chooser as to the correctness of her guess.  (At this point the Guesser's prediction is unknown to the Chooser, but the Chooser does know the amount of the wager.)  Next, the Chooser makes his selection of the node $j$ to move to next.  If the Guesser was correct in her guess, then she receives a payoff proportional to the amount of her wager, while if she was incorrect, she loses the amount of her wager.  The payoff rule is explained in more detail below.  After the payoff both players move to the node $j$ selected by the Chooser.  If node $j$ is a terminal node, then the Guesser's fortune is multiplied by the corresponding value $\val_{j}$, and the game is over.  If node $j$ is non-terminal, play repeats from there.  The goal of the Guesser is to maximize her expected fortune, while the goal of the Chooser is to minimize the Guesser's expected fortune.  We are interested in the optimal strategies for the players in this game, and the dynamics of play under the optimal strategies.

We have adopted a payoff rule which rewards the Guesser for risk-taking.  To motivate this rule, consider first the situation in which the players are at a non-terminal node $i$ of out-degree $\outDeg_{i} \ge 2$.  Then ``on average'' the Guesser will lose her wager $\outDeg_{i}-1$ times for every time that she wins.  This being the case, fairness suggests that her payoff on winning should be $\outDeg_{i}-1$ times her loss on losing.  This \emph{odds-weighted payoff} gives the Guesser incentive to make non-zero wagers at nodes of high out-degree.  Next, consider the case in which the players are at a non-terminal node $i$ of out-degree $\outDeg_{i} = 1$.  The Guesser's prediction will of course be correct in such a case, but the odds-weighted payoff in this case would be zero.  To encourage the Guesser to wager, we choose to suspend odds-weighting in the $\outDeg_{i} = 1$ case, and allow the Guesser to take advantage of a sure bet.  Thus our payoff rule is
\begin{equation} \label{eq:payoffRule}
 \text{new fortune} = 
 \begin{cases}
  \text{current fortune $+ \, (\outDeg_{i}-1) \cdot \text{wager}$} \quad &\text{if $\outDeg_{i} \ge 2$ and Guesser is correct,} \\
  \text{current fortune $+ \,\text{wager}$} \quad &\text{if $\outDeg_{i} = 1$ and Guesser is correct,} \\
  \text{current fortune $- \, \text{wager}$} \quad &\text{if Guesser is incorrect.}
 \end{cases}
\end{equation}
Note that if the out-degree of every vertex is two or less, then the payoff rule is simply: the Guesser wins the amount of her wager if her guess is correct, and loses it otherwise.  

In addition to the payoff rule, optimal strategies obviously depend on the type of graph in which the game is being played.  Call a graph \emph{terminating} if every node has a path to some terminal node.  For terminating graphs, we will show that the values of the non-terminal nodes can be derived from the pre-assigned values of the terminal nodes, and that optimal strategies can be derived from these values.  We will also consider the strongly connected case, which is interesting because there are no terminal nodes, no pre-assigned values, and the game duration is infinite.

%In the next four sections we analyze the Path Guessing Game on fans, trees, terminating graphs, and strongly connected graphs, deriving the optimal strategies and discussing the Markov-chain dynamics of the game along the way.  The last section applies these results to the infinite duration Lying Oracle Game.

Most of our notation is standard.  Random variables and matrices are denoted with uppercase letters; constants and vectors are lowercase.  We will write $i \leadsTo j$ to indicate that there is a directed edge from node $i$ to node $j$ in the graph $\graph$.  The out-degree of vertex $i$ will be denoted by $\outDeg$ or $\outDeg_{i}$.  The symbol $\ones$ will indicate either a vector or matrix, all of whose entries are $1$.  The dimensions should be clear from the contex.  For readability, the Chooser will be consistently referred to as ``he'', and the Guesser as ``she''.

%%%%%%%%%%%%%%%%%%%%%%%%%%%%%%
%%%%%%%%%%%%%%%%%%%%%%%%%%%%%%
%%%%%%%%%%%%%%%%%%%%%%%%%%%%%%
%% fan case
%%%%%%%%%%%%%%%%%%%%%%%%%%%%%%
%%%%%%%%%%%%%%%%%%%%%%%%%%%%%%
%%%%%%%%%%%%%%%%%%%%%%%%%%%%%%
\section{The Game on Fans} \label{s:fans}
Let $\graph$ be a fan with $\outDeg$ leaves (i.e. a tree of height $1$ on $\outDeg+1$ nodes).  The players are initially located at the root node of $\graph$, and each leaf has a positive value $\val_{j}$ assigned to it.
We begin by delineating the strategy sets for each player.

%For consistency with later sections, the root of the fan will be labeled $i$, while the leaves will be indexed by $j$.  
%Each of the leaves of the fan has a positive \emph{value} $\val_{j} > 0$ assigned to it.  The values are public knowledge - both players are aware of them.  
% \begin{center}
%  \includegraphics[width=2.5in]{one-step-game}
% \end{center}
% \begin{center}
%  {Figure 1.  The game on a fan.}
% \end{center}
%\vskip 0.5cm

%The game begins with the Guesser writing down a guess as to which leaf node the Chooser will select.  The guess is kept secret from the Chooser for the time being.  However, the Guesser does announce to the Chooser a wager $\wager$ on the correctness of her guess.  Next, the Chooser chooses a destination node.  The Guesser then reveals her guess, and if she was correct she wins the odds-weighted payoff $\pEnclose{\outDeg-1} \wager$, and her updated fortune is $1 + \pEnclose{\outDeg-1} \wager$.  If the Guesser was incorrect, she loses $\wager$, and her updated fortune is $1 - \wager$.  Finally, her updated fortune is multiplied by the value $\val_{j}$ of the destination node chosen by the Chooser.  The Guesser's goal is to maximize her expected fortune, while the Chooser's goal is to minimize the Guesser's expected fortune.
Let $\cStratSet$ denote the set of pure strategies available to the Chooser.  The Chooser's play consists of selecting the destination node, so
\begin{equation*}
 \cStratSet = \cEnclose{j : 1 \leq j \leq \outDeg}.
\end{equation*}
On the other hand the Guesser must choose both a wager and a guess as to the destination node.  Thus her strategy set $\gStratSet$ is
\begin{equation*}
 \gStratSet = \cEnclose{\pEnclose{j,\wager} : 1 \le j \le  \outDeg \text{ and } \wager \in \bEnclose{0,1}}.
\end{equation*}
Both players are allowed to play non-deterministically.  Let $\cChoice$ denote the (generally random) destination node selected by the Chooser, and similarly let $\gGuess$ and $\gWager$ denote the Guesser's guess and wager respectively.  Naturally we assume that the Chooser's choice $\cChoice$ and Guesser's guess $\gGuess$ are independent: neither the Chooser nor the Guesser are mind readers.  However, the wager is publicly announced, so we leave open the possibility that the Chooser's choice $\cChoice$ is dependent on the Guesser's wager $\gWager$.  Thus, the Chooser's strategy is completely determined by the conditional probabilities
\begin{equation*}
 \cProb(j|\wager) = \cPr{\cChoice = j}{\gWager = \wager}.
\end{equation*}
The Guesser's strategy is specified by the joint distribution of $\gGuess$ and $\gWager$.  This is determined by the conditional distribution of $\gGuess$ given $\gWager$,
\begin{equation*}
 \gProb(j | \wager) = \cPr{\gGuess = j}{\gWager = \wager},
\end{equation*}
along with the unconditional distribution of $\gWager$.

Let $\fort$ denote the Guesser's fortune at the end of the game.  The Chooser's goal is to minimize $\E{\fort}$, while the Guesser wants to maximize $\E{\fort}$.  

\begin{theorem}[Optimal play in the game on fans] \label{thm:optPlayFan}
Let $\graph$ be a fan with $\outDeg$ leaves, each of which has been assigned a positive value $\val_{j}$.  The optimal strategy for the Chooser in the Path Guessing Game on $\graph$ is to make his choice independently of the Guesser's actions, with
\begin{equation} \label{eq:cOptProbFan}
 \cProb(j) = \pr{\cChoice = j} = 
 \begin{cases}
  1 \quad &\text{if $\outDeg = 1$} \\
  \frac{\rValu_{j}}{\sum_{k=1}^{\outDeg} \rValu_{k}} \quad &\text{if $\outDeg \ge 2$}
 \end{cases}
\end{equation}
The optimal strategy for the Guesser is characterized by
\begin{equation} \label{eq:gOptWagerFan0}
 \gWager = 
 \begin{cases}
  1 \quad &\text{if $\outDeg = 1$} \\
  1 - \outDeg \gParam \cProb_{\min} \quad &\text{if $\outDeg \ge 2$}
 \end{cases}
\end{equation}
and
\begin{equation} \label{eq:gOptProbFan0}
 \gProb(j | w) = \cPr{\gGuess = j}{\gWager = \wager} = 
 \begin{cases}
  1 \quad &\text{if $\outDeg = 1$} \\
  \frac{\outDeg \cProb(j) - 1 + \wager}{\outDeg \wager} \quad &\text{if $\outDeg \ge 2$}
 \end{cases}
\end{equation}
where $\cProb_{\min} = \min \cEnclose{\cProb(j) : 1 \le j \le \outDeg}$ and $\gParam$ is any random variable with support in $\ccint{0}{1}$.  Under optimal play, the Guesser's expected fortune at the end of the game is given by
\begin{equation*}
 \E{\fort} = 
 \begin{cases}
  2 \val_{j} \quad &\text{if $\outDeg = 1$} \\
  \frac{\outDeg}{\sum_{j=1}^{\outDeg} \rValu_{j}} \quad &\text{if $\outDeg \ge 2$}
 \end{cases}
\end{equation*}
\end{theorem}

\begin{proof}
The $\outDeg = 1$ case is trivial, so assume $\outDeg \ge 2$.
Let $\hMean$ denote the harmonic mean of the values
\begin{equation} \label{eq:hMeanDef}
 \hMean = \frac{\outDeg}{\sum_{k=1}^{\outDeg} \rValu_{k}}.
\end{equation}
We first show that no matter what strategy the Chooser employs, the Guesser can always find a strategy such that $\E{\fort} \ge \hMean$.  Towards that end, consider any Chooser strategy $\pEnclose{\cProb(k|\wager), k \in \vSet}$.  Then for all $j \in \vSet$ and $\wager \in \ccint{0}{1}$ we have
\begin{align} %\label{eq:cEFortGivenGW}
 \cE{\fort}{\gGuess = j, \gWager = \wager} 
  &= \pEnclose{1 + (\outDeg-1) \wager} \val_{j} \cProb(j | \wager) + \pEnclose{1-\wager} \sum_{k \ne j} \val_{k} \cProb(k | \wager) \notag \\
  &= \outDeg \wager \val_{j} \cProb(j|\wager) + \pEnclose{1-\wager} \sum_{k=1}^{\outDeg} \val_{k} \cProb(k|\wager) \notag \\
  &= \wager \outDeg \val_{j} \cProb(j|\wager) + \pEnclose{1-\wager} \bar{\val} \notag \\
  &= \wager \pEnclose{\outDeg \val_{j} \cProb(j|\wager) - \bar{\val}} + \bar{\val} \notag
\end{align}
where to ease notation $\bar{\val} = \sum_{k=1}^{\outDeg} \val_{k} \cProb(k|\wager)$.
Define the \emph{critical wager}
\begin{equation} \label{eq:critWager0}
 \wagerc = 1 - \hMean \rValu_{\max} \,,
\end{equation}
where $\val_{\max}$ is the maximum value of the nodes in $\graph$.  Fix any $\wager \in \ccint{\wagerc}{1}$, and choose $j = j(\wager)$ such that
\begin{equation} \label{eq:jDef}
 \val_{j} \cProb(j|\wager) \ge \val_{k} \cProb(k|\wager)
\end{equation}
for all $k \in \vSet$.  
We claim that for this $j$ and $\wager$ we have
\begin{equation} \label{eq:cEFortGivenGWEstimate}
 \cE{\fort}{\gGuess = j, \gWager = \wager} \ge \hMean,
\end{equation}
with equality if and only if
\begin{equation} \label{eq:cOptProbFan0}
 \cProb(k | \wager) = \frac{\rValu_{k}}{\sum_{\ell=1}^{\outDeg} \rValu_{\ell}}
\end{equation}
for all $k$, in which case we actually have
\begin{equation} \label{eq:constCEForts}
 \cE{\fort}{\gGuess = k, \gWager = \wager} = \hMean
\end{equation}  
for \emph{all} $k$.  To prove the claim we need to show that
\begin{equation} \label{eq:toShow1}
 \pEnclose{1 - \gParam \hMean \rValu_{\max}} \pEnclose{\outDeg \val_{j} \cProb(j|\wager) - \bar{\val}} + \bar{\val} \ge \hMean
\end{equation}
for all $\gParam \in \ccint{0}{1}$.  But since $j$ satisfies \eqref{eq:jDef}, we have $\outDeg \val_{j} \cProb(j|\wager) - \bar{\val} \ge 0$.  Thus, it is necessary and sufficient to prove \eqref{eq:toShow1} for $\gParam = 1$:
\begin{equation} \label{eq:toShow2}
 \pEnclose{1 - \hMean \rValu_{\max}} \pEnclose{\outDeg \val_{j} \cProb(j|\wager) - \bar{\val}} + \bar{\val} \ge \hMean.
\end{equation}
Let $S = \sum_{k=1}^{\outDeg} \rValu_{k}$, so that $S \hMean = \outDeg$.  Multiplying both sides of \eqref{eq:toShow2} by $S \val_{\max}$ and rearranging, we see that \eqref{eq:toShow1} is equivalent to
\begin{equation} \label{eq:toShow3}
 \cProb(j|\wager) \val_{j} \pEnclose{S \val_{\max} - \outDeg} \ge \val_{\max} - \bar{\val}.
\end{equation}
To establish \eqref{eq:toShow3}, note first that
\begin{equation*}
 S \val_{\max} - \outDeg = \sum_{k=1}^{\outDeg} \frac{\val_{\max}-\val_{k}}{\val_{k}}.
\end{equation*}
Therefore
\begin{align} \label{eq:deriv1}
 \cProb(j|\wager) \val_{j} \pEnclose{S \val_{\max} - \outDeg}
  &= \cProb(j|\wager) \val_{j} \sum_{k=1}^{\outDeg} \frac{\val_{\max}-\val_{k}}{\val_{k}} \notag \\
  &= \sum_{k=1}^{\outDeg} \cProb(k|\wager) \frac{\cProb(j|\wager) \val_{j} \pEnclose{\val_{\max}-\val_{k}}}{\cProb(k|\wager) \val_{k}} \notag \\
  &\ge \sum_{k=1}^{\outDeg} \cProb(k|\wager) \pEnclose{\val_{\max} - \val_{k}} \\
  &= \val_{\max} - \bar{\val} \notag,
\end{align}
where we have used \eqref{eq:jDef} for the inequality.  Therefore we have \eqref{eq:toShow3}, and hence \eqref{eq:cEFortGivenGWEstimate}.  Moreover, there can be equality in \eqref{eq:deriv1} - and hence in \eqref{eq:cEFortGivenGWEstimate} - if and only if
\begin{equation*}
 \cProb(k|\wager) \val_{k} = \cProb(j|\wager) \val_{j}
\end{equation*}
for all $k$, which immediately implies \eqref{eq:cOptProbFan0}.  It is straightforward that \eqref{eq:cOptProbFan0} implies \eqref{eq:constCEForts}, so we have the claim.  
Thus no matter what the Chooser's strategy is, the Guesser can always find a strategy such that $\E{\fort} \ge \hMean$.  Moreover, the claim shows that if the Chooser uses any strategy other than \eqref{eq:cOptProbFan0}, then Guesser can find a strategy such that the inequality is strict.

Now we show that no matter what strategy the Guesser uses, the Chooser can always find a strategy that forces $\E{\fort} \le \hMean$.
Towards that end, note that
\begin{align} \label{eq:cEFortCjWw}
 \cE{\fort}{\cChoice = j, \gWager = \wager} 
  &= \pEnclose{1+(\outDeg-1)\wager} \gProb(j|\wager) \val_{j} + \pEnclose{1-\gProb(j|\wager)} \val_{j} \notag\\
  &= \val_{j} \pEnclose{\wager \pEnclose{\outDeg \gProb(j|\wager) - 1} + 1},
\end{align}
We claim that the expectation in \eqref{eq:cEFortCjWw} is equal to $\hMean$ for all $j$ if and only if 
\begin{equation} \label{eq:derivGWager}
 \wager = 1 - \gParam \hMean \rValu_{\max}
\end{equation}
and
\begin{equation} \label{eq:derivGProb} 
% \gProb(j|\wager) = \frac{\cProb(j) - \gParam \cProb_{\min}}{1 - \outDeg \gParam \cProb_{\min}},
 \gProb(j|\wager) = \frac{\outDeg \cProb(j) - 1 + \wager}{\outDeg \wager},
\end{equation}
where $\cProb(j), j \in \vSet$ are given by \eqref{eq:cOptProbFan}, and $\gParam \in \bEnclose{0,1}$ is arbitrary.
To see this, note that  $\hMean \val^{-1}_{j} = n \cProb(j)$, so the expectation in \eqref{eq:cEFortCjWw} is equal to $\hMean$ for all $j$ if and only if
\begin{equation*}
 \outDeg \cProb(j) - 1 =  \wager \pEnclose{\outDeg \gProb(j|\wager) - 1} \quad \text{for all $j$.}
\end{equation*}
Clearly there is a one-parameter family of solutions to this diagonal system, and direct substitution verifies that \eqref{eq:derivGWager} and \eqref{eq:derivGProb} constitute the solution set.  We further claim that for any other choice of strategy by the Guesser, the Chooser can find a strategy $\cProb$ which makes the Guesser's expected fortune strictly less than $\hMean$.  To see this, note that \eqref{eq:cEFortCjWw} implies that
\begin{equation*}
 \sum_{j=1}^{\outDeg} \frac{\cE{\fort}{\cChoice = j, \gWager = \wager}}{\val_{j}} = \outDeg
\end{equation*}
for any choice of $\wager$ and $\gProb(j | w), j \in \vSet$.  We already know that the expected values in the sum are all equal to $\hMean$ if and only if the Guesser's strategy is given by \eqref{eq:derivGWager} and \eqref{eq:derivGProb}.  Thus, for any other choice of the strategy, some of the expected values in the sum will have to be greater than $\hMean$, and more to the point, some will have to be less than $\hMean$ (so that the value of the sum is still $\outDeg$).  Therefore by concentrating his strategy on the $j$ for which $\cE{\fort}{\cChoice = j, \gWager = \wager} < \hMean$, the Chooser can force $\E{\fort} < \hMean$, as claimed.
This completes the proof.
\end{proof}

A few comments before moving on.  By allowing the players to begin play at \emph{any} node of the fan $\graph$, we can think of the pre-assigned values of the leaf nodes as conditional expectations:
\begin{equation*}
 \val_{j} = \cE{\fort}{\text{Guesser starts at node $j$ with $\$1$}}.
\end{equation*}
Theorem \ref{thm:optPlayFan} then gives the value of the root node $i$:
\begin{equation*}
 \val_{i} = \cE{\fort}{\text{Guesser starts at node $i$ with $\$1$}} = 
 \begin{cases}
  2 &\quad \text{if $\outDeg = 1$} \\
  \frac{\outDeg}{\sum_{j=1}^{\outDeg} \val^{-1}(j)} &\quad \text{if $\outDeg \ge 2$}
 \end{cases}
\end{equation*}
Thus we see that the pre-assigned values of the leaf nodes propagate up the fan to the root node.  This observation will be used repeatedly in the sequel. 

The Chooser's optimal strategy is unique, but Theorem \ref{thm:optPlayFan} shows that there is a family of optimal strategies for the Guesser: the parameter $\gParam$ in \eqref{eq:gOptWagerFan} can have any probability distribution on $\ccint{0}{1}$ whatsoever.  Theorem \ref{thm:optPlayFan} would be simpler if we required that the wager $\gWager$ be deterministic.  Then the optimal strategies for the Guesser would form a one-parameter family, given by
\begin{equation} \label{eq:gOptProbFan}
 \gProb(j) = \pr{\gGuess = j} = \frac{\cProb(j) - \gParam \cProb_{\min}}{1 - \outDeg \gParam \cProb_{\min}}
\end{equation}
and
\begin{equation} \label{eq:gOptWagerFan}
 \wager = 1 - \outDeg \gParam \cProb_{\min}
\end{equation}
The parameter $\gParam \in \ccint{0}{1}$ would then parameterize the Guesser's optimal strategies, and in each one, the wager would be deterministic.  In the interest of simplicity, we will take the wager to be deterministic throughout the rest of the paper.  Accordingly, we drop the random variable notation $\gWager$, and simply denote the wager by $\wager$.  

Note that the critical wager \eqref{eq:critWager0} corresponds to $\gParam = 1$ in \eqref{eq:gOptWagerFan}.  
%We will refer to this as the \emph{critical wager}, and denote it by $\wagerc$:
%\begin{equation} \label{eq:critWager}
% \wagerc = 1 - \outDeg \cProb_{\min} %= 1 - \hMean \rValu_{\max},
%\end{equation}
%%where $\val_{\max}$ is the largest of the leaf node values.  
The corresponding probabilities for the Guesser are
\begin{equation*} \label{eq:minRiskProb}
 \gProb(j) = \frac{\cProb(j) - \cProb_{\min}}{1 - \outDeg \cProb_{\min}}.
\end{equation*}
	We will refer to this as the \emph{minimum risk optimal strategy} for the Guesser.  The \emph{maximum risk optimal strategy} for the Guesser corresponds to $\gParam = 0$ in equations \eqref{eq:gOptProbFan} and \eqref{eq:gOptWagerFan}, in which case the wager is $\wager = 1$ and the Guesser's probabilities are equal to the Chooser's.  Although all the Guesser's strategies are optimal in the sense of maximizing her expected fortune, the maximum risk strategy carries with it the distinct possibility that the Guesser will be bankrupt at the end of the game.  From a practical point of view, the minimum risk strategy is clearly preferable.  Interestingly, in the $\outDeg = 2$ case, if the Guesser follows her minimum risk optimal strategy, then her fortune is guaranteed to equal the harmonic mean of the values: $\fort = \E{\fort} = \hMean$ with probability one.  The proof is left to the interested reader.

%%%%%%%%%%%%%%%%%%%%%%%%%%%%%%
%%%%%%%%%%%%%%%%%%%%%%%%%%%%%%
%%%%%%%%%%%%%%%%%%%%%%%%%%%%%%
%% tree case
%%%%%%%%%%%%%%%%%%%%%%%%%%%%%%
%%%%%%%%%%%%%%%%%%%%%%%%%%%%%%
%%%%%%%%%%%%%%%%%%%%%%%%%%%%%%
\section{The Game on Trees} \label{s:trees}
Let $\graph$ be a tree.  The Path Guessing Game on $\graph$ is a straightforward extension of the game on fans.  Again, the terminal nodes of $\graph$ are the leaves, which we assume have been assigned positive values $\val_{j}$.  The players are initially located at the root.  Note that if play has progressed to the point where the players are at a node $i$ whose children are all leaves, then at that point, they are playing the game on a fan.  The optimal strategies are given by Theorem \ref{thm:optPlayFan}, the only difference being that the Guesser's fortune is some number $\fort$ rather than one dollar.  But we can account for this by thinking of the wager $\wager \in \bEnclose{0,1}$ as the \emph{proportion} of the Guesser's fortune that is wagered.  If the out-degree of node $i$ is at least $2$, then by Theorem \ref{thm:optPlayFan} the Guesser's expected fortune at the end of the game, given that her fortune at node $i$ is $\fort$, is $\hMean_{i} \fort$, where $\hMean_{i}$ is the harmonic mean of the values of the leaves of node $i$.  Thus we should assign the value of node $i$ to be $\val_{i} =\hMean_{i}$ in this case.  In the same way, if the out-degree of $i$ is $1$, then Theorem \ref{thm:optPlayFan} implies $\val_{i} = 2 \val_{\ell}$, where $\ell$ is the single child of $i$. Likewise, values can be assigned to all the interior nodes of the tree, all the way up to the root.  Playing in accordance with these values is optimal, by a straightforward induction on the height of the tree.  These ideas are summarized in the following theorem.

\begin{theorem}[Optimal play on trees]\label{thm:optPlayTree}
Let $\graph$ be a finite rooted tree in which each leaf node $\leafNode$ has been assigned a positive value $\val_{\leafNode}$.  Define the values of all the other nodes $i$ of $\graph$ by
\begin{equation} \label{eq:treeValues}
 \rValu_{i} = 
 \begin{cases}
  \rValu_{j} / 2 \quad &\text{if $\outDeg_i = 1$ and $i \leadsTo j$} \\
  \outDeg_{i}^{-1} \sum_{j : i \leadsTo j} \rValu_{j} \quad &\text{if $\outDeg_i \ge 2$}
 \end{cases}
\end{equation}
The optimal strategy for the Chooser is given by
\begin{equation} 
 \cProb_{i,j} = 
 \begin{cases}
  1 \quad &\text{if $\outDeg_{i} = 1$ and $i \leadsTo j$} \\
  \rValu_{j} / \sum_{k : i \leadsTo k} \rValu_{k} \quad &\text{if $\outDeg_{i} \ge 2$ and $i \leadsTo j$} \\
  0 \quad &\text{otherwise}
 \end{cases}
\end{equation}
There is a one-parameter family of optimal strategies for the Guesser, given by
\begin{equation} 
 \gProb_{i,j} = 
 \begin{cases}
  1 \quad &\text{if $\outDeg_{i} = 1$ and $i \leadsTo j$} \\
  \pEnclose{\cProb_{i,j} - \gParam \cProb_{i,\min}}/ \pEnclose{1 - \outDeg_{i} \gParam \cProb_{i,j}} \quad &\text{if $\outDeg_{i} \ge 2$ and $i \leadsTo j$} \\
  0 \quad &\text{otherwise}
 \end{cases}
\end{equation}
and
\begin{equation}
 \wager_{i} = 
 \begin{cases}
  1 & \quad \text{if $\outDeg_{i} = 1$} \\
  1 - \outDeg_{i} \gParam \cProb_{i,\min} & \quad \text{if $\outDeg_{i} \ge 2$} 
 \end{cases}
\end{equation}
where $\cProb_{i,\min} = \min \cEnclose{\cProb_{i,j} : \text{$j$ such that $i \leadsTo j$}}$, and $\gParam$ is arbitrary in $\bEnclose{0,1}$.  Under optimal play the Guesser's fortune $\fort$ at the end of the game satisfies $\E{\fort} = \val_{0}$, where $\val_{0}$ is the value of the root node of the tree.
\end{theorem}

%\begin{example}
%As an example, consider the tree in Figure 2 (\textit{a}), in which the leaf nodes have been assigned values of $2$, $1$, and $1$.  Moving up one level in the tree in part (\textit{b}) of the figure, the node whose leaves have values of $2$ and $1$ must itself have a value of $4/3$, the harmonic mean of $2$ and $1$.  The node whose single child has value of $1$ will have a value of $2$, since the Guesser is guaranteed to double her fortune from that point forward.  Finally, we can calculate the value of the root node of the tree: it is $8/5$, the harmonic mean of $4/3$ and $2$, as shown in part (\textit{c}).

% \begin{center}
%  \includegraphics[width=5in]{value-propagates}
% \end{center}
%  Figure 2.  Value propagates.  In (\textit{a}) the leaf nodes are arbitrarily assigned values of $2$, $1$, and $1$.  In (\textit{b}) and (\textit{c}) the values are propagated up the tree to the root node.  The total value of the game to the Guesser is $1.6$, the value of the root node.  The optimal strategies for the game on the tree of Figure 2 are shown in Figure 3.
% \begin{center}
%  \includegraphics[width=5in]{opt-strats-tree-ex}
% \end{center}
% Figure 3.  Optimal strategies for (\textit{a}) the Chooser and (\textit{b}) the Guesser.  The parameter $\gParam$ is arbitrary in $\bEnclose{0,1}$.
% \vskip 0.5cm
%\hfill \expSymb
%\end{example}

%%%%%%%%%%%%%%%%%%%%%%%%%%%%%%
%%%%%%%%%%%%%%%%%%%%%%%%%%%%%%
%%%%%%%%%%%%%%%%%%%%%%%%%%%%%%
%% terminating case
%%%%%%%%%%%%%%%%%%%%%%%%%%%%%%
%%%%%%%%%%%%%%%%%%%%%%%%%%%%%%
%%%%%%%%%%%%%%%%%%%%%%%%%%%%%%
\section{The Game on Terminating Graphs} \label{s:terminating}
A connected digraph in which every non-terminal node has a (directed) path to some terminal node will be referred to as a \emph{terminating graph}.  
In this section we assume that $\graph$ is a terminating graph with $\nNodes$ vertices, in which each terminal node $j$ has been pre-assigned a positive value $\val_{j}$. Note that since a terminating graph may contain cycles and other strongly connected subgraphs, the duration of play is potentially infinite.  We will derive optimal play in the Path Guessing Game on $\graph$ by first assigning appropriate values to the non-terminal nodes, then proving that playing the game in accordance with these values (as in Theorem \ref{thm:optPlayFan}) is optimal.

We begin by considering a truncated game which lasts for at most $\tmNds$ steps.  The truncated game is identical to the non-truncated game, except that if the players have not reached a terminal node by step $\tmNds$, the game is stopped.  We will refer to this as the \emph{$\tmNds$-step game} (even though it may be over in fewer than $\tmNds$ steps).  The $\tmNds$-step game can be completely described by $\nNodes_{\ntl}$ \emph{path trees}, where $\nNodes_{\ntl}$ is the number of non-terminal nodes of $\graph$.  The path tree corresponding to a non-terminal vertex $i$ represents all the paths emanating from $i$ whose lengths are at most $\tmNds$.  Terminal nodes in the original graph $\graph$ will be leaf nodes in the path trees.  The non-terminal nodes that are reachable in $\tmNds$ steps from node $i$ will also appear as leaf nodes in the path tree corresponding to $i$.  For the path trees in this truncated game, terminal nodes retain the values assigned to them in the non-truncated game.  Leaf nodes of the path trees that correspond to non-terminal nodes are assigned the value of $1$, in keeping with the idea that the game simply stops if one of these nodes is reached.  
%Figure 4 below illustrates the idea.
% \begin{center}
%  \includegraphics[width=5.5in]{path-trees-1}
% \end{center}
%  Figure 4.  (\textit{a}) A terminating graph, in which the terminal node has a value of $2$.  (\textit{b}) the corresponding path trees of length $\tmNds = 3$.  Leaves that are not terminal nodes are assigned a value of $1$.
% \vskip 0.5cm 
Thus we have reduced the $\tmNds$-step game to a game on the path trees of length $\tmNds$, which can be solved by propagating the values up the path trees using Theorem \ref{thm:optPlayTree}.

The propagation of the values for the $\tmNds$-step game can be summarized nicely using a matrix approach.  Let $\nNodes_{\ntl}$ be the number of non-terminal nodes, and let $\nNodes_{\tl}$ be the number of terminal nodes.  Number the nodes so that the non-terminal nodes are numbered first, followed by the terminal nodes.  Consider first the path trees for a $1$-step game on $\graph$.  
%In each path tree, the leaves that are terminal nodes of $\graph$ retain their pre-assigned values, while the leaves that are non-terminal nodes of $\graph$ are assigned a value of $1$.
Let $\rVal_{0}$ denote the vector of the reciprocal values of the leaf nodes.  In keeping with our numbering convention, $\rVal_{0}$ is a partitioned vector of the form
\begin{equation} \label{eq:leafNodeVals}
 \rVal_{0} =
 \left( 
  \begin{array}{c}
   \ones \\ \hline
   \rVal_{\tl}
  \end{array}
 \right),
\end{equation}
where $\ones$ is the $\nNodes_{\ntl} \times 1$ vector of ones, and $\rVal_{\tl}$ is the $\nNodes_{\tl} \times 1$ vector containing the pre-assigned reciprocal values of the terminal nodes of $\graph$.
By Theorem \ref{thm:optPlayFan} the reciprocal values of the $\nNodes_{\ntl}$ root nodes of the path trees for the $1$-step game are given by
\begin{equation*} \label{eq:valProp1}
 \rVal_{1} = \propM \rVal_{0},
\end{equation*}
where $\propM$ is the $\nNodes \times \nNodes$ \emph{propagation matrix} $\propM = (\propE_{i,j}: i, j \in \vSet)$ defined by
\begin{equation} \label{eq:propMDef}
 \propE_{i,j} = 
 \begin{cases}
  1               \quad &\text{if $i$ is a terminal node and $i = j$} \\
  1 / 2           \quad &\text{if $\outDeg_{i} = 1$ and $i \leadsTo j$} \\
  1 / \outDeg_{i} \quad &\text{if $\outDeg_{i} \ge 2$ and $i \leadsTo j$} \\
  0               \quad &\text{otherwise}
 \end{cases}.
\end{equation}
%Next consider the path trees of length two, corresponding to playing the $2$-step game on $\graph$.  The values of the leaf nodes are again given by \eqref{eq:leafNodeVals}.  The values of the nodes at one level above the leaf nodes are given by \eqref{eq:valProp1}, because starting from such a node, we are playing the game for one step on $\graph$.  Thus, the reciprocal values of the root nodes of the path trees of length two are given by
%\begin{equation*} \label{eq:valProp2}
% \rVal_{2} = \propM \rVal_{1} = \propM^{2} \rVal_{0}.
%\end{equation*}
Inductively, the reciprocal values for the root nodes of the path trees for the $\tmNds$-step game are given by
\begin{equation} \label{eq:valProp3}
 \rVal_{\tmNds} = \propM^{\tmNds} \rVal_{0}.
\end{equation}
Our first goal is to prove that $\lim_{\tmNds \approaches \infty} \rVal_{\tmNds}$ exists.  In accordance with our numbering convention, the partitioned form of the propagation matrix $\propM$ is
\begin{equation} \label{eq:propMPartition}
 \propM = 
 \left[ 
  \begin{array}{c|c}
   \propMntnt & \propMntt \\ \hline
   0 & \id
  \end{array}
 \right],
\end{equation}
where $\propMntnt$ is $\nNodes_{\ntl} \times \nNodes_{\ntl}$, $\propMntt$ is $\nNodes_{\ntl} \times \nNodes_{\tl}$, the zero submatrix is $\nNodes_{\tl} \times \nNodes_{\ntl}$, and the identity submatrix is $\nNodes_{\tl} \times \nNodes_{\tl}$.  It follows that the powers of $\propM$ are of the form
\begin{equation} \label{eq:propMPowers}
 \propM^{\tmNds} = 
 \left[ 
  \begin{array}{c|c}
   \propMntnt^{\tmNds} & (\sum_{i=0}^{\tmNds-1} \propMntnt^{i}) \propMntt \\ \hline
   0 & \id
  \end{array}
 \right].
\end{equation} 
The next lemma establishes that limiting values exist, are strictly positive, and that the limiting reciprocal value vector $\rVal $ satisfies $\propM \rVal = \rVal$.
\begin{lemma}[Limiting Values] \label{lem:limitingValuesTerm}
 Let $\graph$ be a terminating graph, and let $\propM$ be the associated propagation matrix defined by \eqref{eq:propMDef}.  Then $\lim_{\tmNds \approaches \infty} \propM^{\tmNds}$ exists, with
 \begin{equation} \label{eq:propMLimit}
  \lim_{\tmNds \approaches \infty} \propM^{\tmNds} =
  \left[ 
   \begin{array}{c|c}
    0 & (\id - \propMntnt)^{-1} \propMntt \\ \hline
    0 & \id
   \end{array}
  \right],
 \end{equation}
where $\propMntnt$ and $\propMntt$ are the submatrices defined by \eqref{eq:propMPartition}.  In particular, the limiting reciprocal value vector $\rVal = \lim_{\tmNds \approaches \infty} \propM^{\tmNds} \rVal_{0}$ exists, satisfies $\propM \rVal = \rVal$, and is strictly positive.  Partitioning the limiting reciprocal value vector consistently with the partition of $\propM$,
\begin{equation} \label{eq:rValLimit}
 \rVal =
 \left( 
  \begin{array}{c}
   \rVal_{\ntl} \\ \hline
   \rVal_{\tl}
  \end{array}
 \right),
\end{equation}
the limiting reciprocal values $\rVal_{\ntl}$ of the non-terminal nodes are related to the values $\rVal_{\tl}$ of the terminal nodes by
\begin{equation} \label{eq:rValntl}
 \rVal_{\ntl} = \pEnclose{\id - \propMntnt}^{-1} \propMntt \rVal_{\tl},
\end{equation}
\end{lemma}
\begin{proof}
Note that if all the eigenvalues of $\propMntnt$ were less than $1$ in absolute value, we could immediately conclude that $\lim_{\tmNds \approaches \infty} \sum_{i=0}^{\tmNds-1} \propMntnt^{i} = (\id - \propMntnt)^{-1}$, and the result would follow.  We claim that in fact all the eigenvalues of $\propMntnt$ \emph{are} strictly less than $1$ in absolute value.  To see this, first observe that $\propMntnt$ is a nonnegative sub-stochastic matrix, so by the standard Perron-Frobenius theory of nonnegative matrices (see \cite{Minc}, for instance), the maximal eigenvalue $\propMEV$ of $\propMntnt$ is nonnegative and less than or equal to $1$, and $\abs{\lambda} \le \propMEV$ for any other eigenvalue $\lambda$ of $\propMntnt$.  My claim is that $r < 1$.  Suppose (by way of contradiction) that $r = 1$.  Then there is an associated nonnegative eigenvector $x$, which we may assume without loss of generality has been scaled so that its maximum component is equal to $1$.  Now the equation $\propMntnt x = x$ along with the definition \eqref{eq:propMDef} imply the following:
\begin{enumerate}
 \item \label{it:xItem1} If the non-terminal node $i$ of $\graph$ has out-degree $1$ (i.e. $\outDeg_{i} = 1$), and if the (sole) child of $i$ is a terminal node $j$, then $x_{i} = 0$. 
 \item \label{it:xItem2} If the non-terminal node $i$ has $\outDeg_{i} = 1$, and if its (sole) child is a non-terminal node $j$, then $x_{i} = x_{j} / 2$.
 \item \label{it:xItem3} If the non-terminal node $i$ has $\outDeg_{i} \ge 2$, then 
  \begin{equation} \label{eq:xItem3}
   x_{i} = \outDeg_{i}^{-1} \sum_{j \in J(i)} x_{j},
  \end{equation}
  where $J(i) = \cEnclose{\text{non-terminal nodes $j$ such that $\propE_{i,j} \ne 0$}}$.
\end{enumerate}
Items \ref{it:xItem1} and \ref{it:xItem2} imply that entries of $x$ corresponding to non-terminal nodes of out-degree $1$ are all strictly less than $1$.  Now consider the remaining entries of $x$, which correspond to non-terminal nodes $i$ of out-degree $2$ or more.  Let $\dist_{i}$ denote the distance from non-terminal node $i$ to the nearest terminal node.  Item \ref{it:xItem1} says that if $\dist_{i} = 1$, then $x_{i} = 0$.  Now suppose $\dist_{i} = 2$.  Then there is a directed edge from node $i$ to some  node $j$ with $\dist_{j} = 1$.  Thus $x_{j} < 1$, and therefore $x_{i} < 1$ by \eqref{eq:xItem3} (since the cardinality of $J(i)$ is at most $\outDeg_{i}$).  So any non-terminal node $i$ with $\dist_{i} = 2$ must have $x_{i} < 1$.  An easy induction shows that any non-terminal node $i$ that is at a finite distance from some terminal node must have $x_{i} < 1$.  But in a terminating graph \emph{every} non-terminal node is at a finite distance from some terminal node.  Therefore $x_{i} < 1$ for \emph{all} non-terminal nodes $i$, contradicting our scaling of $x$ so that its maximum component was $1$.  Thus $r < 1$ as claimed, and we have equation \eqref{eq:propMLimit}.

From \eqref{eq:valProp3} and \eqref{eq:propMLimit} it follows that the limiting reciprocal value $\rVal = \lim_{\tmNds \approaches \infty} \rVal_{\tmNds}$ exists, and satisfies
\begin{equation} \label{eq:rValLimit0}
 \rVal = 
 \left[ 
  \begin{array}{c|c}
   0 & (\id - \propMntnt)^{-1} \propMntt \\ \hline
   0 & \id
  \end{array}
 \right]
 \rVal_{0}.
\end{equation}
Equation \eqref{eq:rValntl} follows immediately.  Using \eqref{eq:rValntl} and \eqref{eq:propMPartition} it is easy to see that the limiting vector $\rVal$ is a right eigenvector of the propagation matrix, corresponding to the maximal eigenvalue $\propMEV = 1$.
%\begin{align} \label{eq:rValEValTerm}
%	\propM \rVal = \rVal.
%\end{align}

It remains to show that $\rVal > 0$.  We claim that no row of $\pEnclose{\id - \propMntnt}^{-1} \propMntt$ consists entirely of zeros.  To see this, note that each non-terminal node $i$ is connected to some terminal node $k$ by a path in $\graph$ of some length $\tmNds > 0$.  The node $j$ immediately preceeding $k$ in this path must be a non-terminal node.  It follows that the $i,j$-entry in $\propMntnt^{s-1}$ is nonzero, as is the $j,k$-entry in $\propMntt$.  Therefore row $i$ of $A^{s-1} B$ can not consist entirely of zeros.  Since $\pEnclose{\id - \propMntnt}^{-1} \propMntt$ is a sum of such (nonnegative) terms, its $\ith{i}$ row can not consist entirely of zeros either, establishing the claim.  Now the pre-assigned reciprocal value vector $\rVal_{\tl}$ is positive by assumption, so by  \eqref{eq:rValntl} we conclude that $\rVal_{\ntl}$ is strictly positive.  This completes the proof of the lemma.
\end{proof}

Now we turn to the full, non-truncated game.  The previous lemma shows that for a terminating graph, any assignment of values for the terminal nodes leads to a unique set of limiting values for the non-terminal nodes.  A reasonable option for the players in the non-truncated game is to play in accordance with these limiting values.  For the Chooser, this would mean
\begin{equation} \label{eq:cOptProbTerm}
 \cProb_{i,j} = \cPr{\pos_{\tmNdx+1} = j}{\pos_{\tmNdx} = i} = 
 \begin{cases}
  1 &\quad \text{if $\outDeg_i = 1$} \\
  \rVal_{j} / \pEnclose{\sum_{k : i \leadsTo k} \rVal_{k}} &\quad \text{if $\outDeg_i \ge 2$}
 \end{cases},
\end{equation}
where $\pos_{\tmNdx}$ denotes the vertex occupied by the players at time $\tmNdx$.  For the Guesser, it would mean
\begin{equation} \label{eq:gOptProbTerm}
 \gProb_{i,j} = \cPr{\gGuess_{\tmNdx+1} = j}{\pos_{\tmNdx} = i} = 
 \begin{cases}
  1 &\quad \text{if $\outDeg_i = 1$} \\
  \pEnclose{\cProb_{i,j} - \gParam \cProb_{i,\min}} / \pEnclose{1 - \outDeg_{i} \gParam \cProb_{i,\min}} &\quad \text{if $\outDeg_i \ge 2$}
 \end{cases}
\end{equation}
where $\gGuess_{\tmNdx+1}$ is the Guesser's prediction of where the Chooser will go on the next step, with the wager at node $i$ given by
\begin{equation} \label{eq:gOptWagerTerm}
 \wager_{i} = 
 \begin{cases}
  1 & \quad \text{if $\outDeg_{i} = 1$} \\
  1 - \outDeg_{i} \gParam \cProb_{i,\min} & \quad \text{if $\outDeg_{i} \ge 2$} 
 \end{cases},
% \cPr{\gWager_{\tmNdx+1} = 1 - \outDeg_{i} \gParam \cProb_{i,\min}}{\pos_{\tmNdx} = i} = 1 
\end{equation}
where $\gParam$ is arbitrary in $\ocint{0}{1}$, and $\cProb_{i,\min}$ is the probability of the vertex least likely to be chosen by the Chooser:
\begin{equation*}
 \cProb_{i,\min} = \min \cEnclose{\cProb_{i,j} : \text{$j \in \vSet$, $i \leadsTo j$}}
\end{equation*}
%As in Section \ref{s:fans}, We have used $\wager_{i}$ in place of $\gWager_{i}$, since the optimal wager at node $i$ is deterministic.  (Recall also that $\wager_{i}$ represents the \emph{fraction} of the Guesser's current holdings that she wagers.)  
We will refer to the strategies given by \eqref{eq:cOptProbTerm}, \eqref{eq:gOptProbTerm}, and \eqref{eq:gOptWagerTerm} as the \emph{limiting strategies} for the game.

Let $\cProbM = \pEnclose{\cProb_{i,j}: i,j \in \vSet}$ denote the transition probability matrix corresponding to the Chooser's limiting strategy \eqref{eq:cOptProbTerm}.  These are the transition probabilities for the Markov chain that describes the players' position in the graph as a function of time.  Strictly speaking, these probabilities only make sense if $i$ is a non-terminal node.  However, we can think of each terminal node as having an attached loop leading back to itself, with the understanding that all wagering and guessing stops as soon as the players reach a terminal node.  Then $\cProb_{i,i} = 1$ if $i$ is a terminal node, which makes the Markov chain $\pos_{\tmNdx}$ well-defined for all $\tmNdx$.  The end of actual play in the game is given by the \emph{hitting time} $\sTime$ of the chain to the set of terminal nodes:
\begin{equation*} %\label{eq:sTimeDef}
 \sTime = \min \cEnclose{\tmNdx \ge 0 : \text{$\pos_{\tmNdx}$ is a terminal node}}.
\end{equation*}
The next theorem gives the basic properties of the game under the limiting strategies.
\begin{theorem}[Limiting strategies] \label{thm:limitingStrategiesTerm}
Let $\graph = \pEnclose{\vSet, \eSet}$ be a terminating graph in which each terminal node has been assigned a positive value.  Let $\propM$ be the propagation matrix for $\graph$ defined by \eqref{eq:propMDef}, and let $\val = \pEnclose{\val_{i} : i \in \vSet}$ be the vector of limiting values.  Then under the limiting strategies the Chooser's transition probability matrix $\cProbM$ is diagonally similar to the propagation matrix:
\begin{equation} \label{eq:cOptProbMTerm}
 \cProbM = \valM \propM \valM^{-1},
\end{equation}
where $\valM = \text{diag}(\val)$ is the diagonal matrix of limiting values. Under the limiting strategies, the game duration $\sTime$ is finite with probability $1$, and the Guesser's expected fortune at the end of the game is equal to the value of the node at which play began:
\begin{equation} \label{eq:eFortTerm}
 \cE{\fort_{\sTime}}{\pos_{0} = i} = \val_{i}.
\end{equation} 
\end{theorem}
\begin{proof}
Because the  vector $\rVal$ of limiting reciprocal values satisfies $\propM \rVal = \rVal$, equation \eqref{eq:cOptProbTerm} reduces to
\begin{equation*} 
\cProb_{i,j} = \propE_{i,j} \, \rVal_{j} / \rVal_{i} = \propE_{i,j} \, \val_{i} / \val_{j}
\end{equation*}
for \emph{all} $i$ and $j$, and \eqref{eq:cOptProbMTerm} follows.
Because the limiting reciprocal value vector $\rVal$ is strictly positive, it follows from \eqref{eq:cOptProbTerm} that $\cProb_{i,j} > 0$ whenever $i \leadsTo j$ in $\graph$.  Therefore, the hitting time $\sTime$ is finite with probability $1$ under the limiting strategies.

For \eqref{eq:eFortTerm}, note first that it is obviously true if $i$ is a terminal node.  If $i$ is non-terminal, then conditioning on the position of the players at time $\tmNdx = 1$ gives
\begin{equation} \label{eq:eFortTerm1}
 \cE{\fort_{\sTime}}{\pos_{0} = i} = 
 \begin{cases}
  2 \cE{\fort_{\sTime}}{\pos_{0} = j}  & \quad \text{if $\outDeg_{i} = 1$ and $i \leadsTo j$} \\
  \sum_{j:i \leadsTo j} \cProb_{i,j} \cE{\fort_{\sTime}}{\pos_{0} = i,\pos_{1} = j} & \quad \text{if $\outDeg_{i} \ge 2$}
 \end{cases}.
\end{equation}
Now in the second case we have
\begin{align*}
 \cE{\fort_{\sTime}}{\pos_{0} 
  = i,\pos_{1} = j} 
  &= \bEnclose{ \gProb_{i,j} \pEnclose{1+(\outDeg_{i}-1) \wager_{i}} + \pEnclose{1-\gProb_{i,j}} \pEnclose{1 - \wager_{i}} } \cE{\fort_{\sTime}}{\pos_{0} = j} \\
  &= \pEnclose{\wager_{i} \pEnclose{\outDeg_{i} \gProb_{i,j}-1}+1} \cE{\fort_{\sTime}}{\pos_{0} = j} \\
  &= \outDeg_{i} \cProb_{i,j} \cE{\fort_{\sTime}}{\pos_{0} = j} \\
  &= \frac{\val_{i}}{\val_{j}} \cE{\fort_{\sTime}}{\pos_{0} = j},
\end{align*}
so that \eqref{eq:eFortTerm1} becomes
\begin{equation} \label{eq:eFortTerm2}
 \cE{\fort_{\sTime}}{\pos_{0} = i} =
 \begin{cases}
  2 \cE{\fort_{\sTime}}{\pos_{0} = j}  & \quad \text{if $\outDeg_{i} = 1$ and $i \leadsTo j$} \\
  \sum_{j:i \leadsTo j} \propE_{i,j} \frac{\val_{i}^{2}}{\val_{j}^{2}} \cE{\fort_{\sTime}}{\pos_{0} = j} & \quad \text{if $\outDeg_{i} \ge 2$}
 \end{cases}.
\end{equation}
This is recognized as the matrix equation
\begin{equation*}
 \zeta = \valM^{2} \propM \valM^{-2} \zeta,
\end{equation*}
where $\zeta(i) = \cE{\fort_{\sTime}}{\pos_{0} = i}$.  So the vector $z = \valM^{-2} \zeta$ satisfies $\propM z = z$, and using the partition in \eqref{eq:propMPartition} we get
\begin{equation} \label{eq:zNtTerm}
 \propMntnt z_{\ntl} + \propMntt z_{\tl} = z_{\ntl},
\end{equation}
where we have partitioned $z$ consistent with the partition of $\propM$ in \eqref{eq:propMPartition}.  We know that for \emph{terminal} nodes $i$ we have $\zeta(i) = \val_{i}$, and it follows that $z_{\tl} = \valM_{\tl}^{-2} \zeta_{\tl} = \rVal_{\tl}$.  Therefore the solution $z_{\ntl}$ of \eqref{eq:zNtTerm} is $z_{\ntl} = \pEnclose{\id - \propMntnt}^{-1} \propMntt \rVal_{\tl}$, which equals $\rVal_{\ntl}$ by \eqref{eq:rValntl}.  Thus $z = \rVal$, and $\zeta = \valM^{2} z = \val$.  This completes the proof.
\end{proof}

Other quantities of interest in the game include the terminal probabilities
\begin{equation*}
 \tProb_{i,k} = \cPr{\pos_{\sTime} = k}{\pos_{0} = i},
\end{equation*}
the expected stopping times
\begin{equation*} %\label{eq:eSTimeDef}
 \eSTime_{i} = \cE{\sTime}{\pos_{0} = i},
\end{equation*}
and the conditional distributions of the stopping time
\begin{equation*} %\label{eq:sTimeProbDef}
 \sTimeProb_{\tmNdx | i} = \cPr{\sTime = \tmNdx}{\pos_{0} = i}.
\end{equation*}
Form the matrix 
$$\tProb = \pEnclose{\tProb_{i,k}: \text{$i$ is non-terminal, $j$ is terminal}},$$
and the vectors 
$$\sTimeProb_{\tmNdx} = \pEnclose{\sTimeProb_{\tmNdx | i} : \text{$i$ is non-terminal}}$$ 
and 
$$\eSTime = \pEnclose{\eSTime_{i} : \text{$i$ is non-terminal}}.$$  
(In each case the values when $i$ is terminal or $j$ is non-terminal are obvious.)  The next theorem expresses these quantities in terms of the propagation matrix and the limiting values.  In the statement of the theorem, $\valM_{\ntl} = \text{diag}(\val_{\ntl})$ and  $\valM_{\tl} = \text{diag}(\val_{\tl})$ are the diagonal matrices of limiting non-terminal and terminal values respectively.
\begin{theorem} \label{thm:paramsOfInterestTerm}
Let $\graph$ be a terminating graph, and suppose the players use the limiting strategies defined by \eqref{eq:cOptProbTerm}, \eqref{eq:gOptProbTerm}, and \eqref{eq:gOptWagerTerm}.  Then the vector $\eSTime$ of expected stopping times is given by
\begin{equation*} %\label{eq:sTimeDist}
 \eSTime = \valM_{\ntl} \pEnclose{\id - \propMntnt}^{-1} \valM^{-1}_{\ntl} \ones.
\end{equation*}
The conditional stopping time distributions are given by
\begin{equation*} %\label{eq:sTimeProb}
 \sTimeProb_{\tmNdx} = \valM_{\ntl} \propMntnt^{\tmNdx - 1} \propMntt \valM^{-1}_{\tl} \ones,
\end{equation*}
and the terminal probabilities are
\begin{equation*} %\label{eq:tProb}
 \tProb = \valM_{\ntl} \pEnclose{\id - A}^{-1} \propMntt \valM^{-1}_{\tl}.
\end{equation*}
\end{theorem}
\begin{proof}
The proof proceeds by conditioning on the player's position at time $\tmNdx = 1$.  Details are left to the reader.
\end{proof}

As might be expected, the limiting strategies are in fact optimal for the game on terminating graphs.

\begin{theorem}[Optimal play on terminating graphs] \label{thm:optPlayTerm}
Let $\graph$ be a terminating graph in which each terminal node has been assigned a positive value.  Then the limiting strategies \eqref{eq:cOptProbTerm}, \eqref{eq:gOptProbTerm}, and \eqref{eq:gOptWagerTerm} are optimal.
\end{theorem}

\begin{proof}
It is clear that we can restrict our attention to strategies that are \emph{purely positional}, in the sense that at every time the players are at a given vertex $i$, they play the same strategy.  
We will continue to use $E$ to denote expected values under the limiting strategies, while $E^{*}$ will denote expected values under general (but fixed) purely positional strategies.  Suppose both players are playing a general homogeneous strategy, possibly different from the limiting strategies defined by \eqref{eq:cOptProbTerm}, \eqref{eq:gOptProbTerm}, and \eqref{eq:gOptWagerTerm}.  
Define $\val^{*}(i)$ by
\begin{equation} \label{eq:valueLimSup}
 \val^{*}(i) = \limsup_{\tmNdx \approaches \infty}  E^{*} \bEnclose{\fort_{\tmNdx} | \pos_{0} = i}.
\end{equation}
Note that for terminal vertices $i$, $\val^{*}(i)$ equals the pre-assigned value $\val_{i}$. We first claim that if the players are playing optimally, then $0 < \val^{*}(i) < \infty$ for all non-terminal vertices $i$.  Indeed, $\val^{*}(i) > 0$ because the Guesser can always elect to bet zero at every vertex, while $\val^{*}(i) < \infty$ because the Chooser can always elect to take the shortest route from vertex $i$ to a terminal vertex. 

%%%%%%%%%%%%%%%%%%%%%%%%%%%%%%
%% detailed proof the T is finite with probability 1
%%%%%%%%%%%%%%%%%%%%%%%%%%%%%%
Next we claim that if the Chooser is playing optimally, then $\sTime$, the hitting time to the set of terminal states, must be finite with probability 1.
To see this, suppose by way of contradiction that $P^{*}(\sTime = \infty) > 0$.  Then there must exist a non-terminal vertex $i$ that is part of a strongly connected subgraph $\graph^{*}$ of $\graph$ that can be visited infinitely often by the players.  Without loss of generality there is an edge from $i$ to a vertex $j$ not in $\graph^{*}$ that leads to some terminal node $z$, and such that the Chooser's probability of selecting the edge from $i$ to $j$ is zero.  (Otherwise, $\sTime$ would be finite with probability $1$.)  
% \begin{center}
%  \includegraphics[width=3in]{T-finite-fig}
% \end{center}
% \begin{center}
%  Figure 5.  
% \end{center}
The out-degree $\outDeg_{i}$ of $i$ is therefore at least $2$.  Consider the following strategy for the Guesser: when at vertex $i$, the Guesser bets one half of her fortune on the node $k$ that the Chooser is most likely to select; when at any other vertex in $\graph^{*}$, the Guesser bets nothing.  Because the Chooser goes from $i$ to $j$ with probability $0$, the probability $\cProb^{*}_{i,k}$ that he goes to the node $k$ must satisfy
$$\cProb^{*}_{i,k} \ge \frac{1}{\outDeg_{i} - 1}.$$
Therefore the Guesser's expected fortune after playing at vertex $i$ (as a proportion of her current fortune) are
\begin{align*}
 \frac{1}{2} \pEnclose{1 - \cProb^{*}_{i,k}} + \pEnclose{1 + \frac{1}{2} \pEnclose{\outDeg_{i}-1}} \cProb^{*}_{i,k}
 &= \frac{1}{2} \pEnclose{\cProb^{*}_{i,k} \outDeg_{i} + 1} \\
 &\ge \frac{1}{2} \pEnclose{\frac{\outDeg_{i}}{\outDeg_{i}-1} + 1} \\
 &> 1.
\end{align*}
Thus the Guesser increases her fortune on average each time the players are at vertex $i$.  Since there is a positive probability this will happen infinitely often, $\val^{*}(i)$ as defined by \eqref{eq:valueLimSup} will be infinite.  But this contradicts the claim proven above that $\val^{*}(i)$ is finite if the Chooser is playing optimally.  Therefore the claim that $P^{*}(\sTime = \infty) = 0$ is established.

Now $\sTime$ is finite with probability $1$, it follows that $\val^{*}(i)$ as defined by \eqref{eq:valueLimSup} satisfies
\begin{equation*}
 \val^{*}(i) = E^{*} \bEnclose{\fort_{\sTime} | \pos_{0} = i}.
\end{equation*}
provided that the players are playing optimally.
Therefore, given that play starts at node $i$, the game on $\graph$ is equivalent to the game on the fan whose leaves have the values $\val^{*}(j)$, where $i \leadsTo j$.
% \begin{center}
%  \includegraphics[width=2in]{equiv-one-step-term}
% \end{center}
% \begin{center}
%  Figure 6.  
% \end{center}
Optimal play in this game is given by Theorem \ref{thm:optPlayFan}.  Since this is true for each vertex $i$ in $\graph$, the reciprocal value vector $\rVal^{*} = \pEnclose{\rVal^{*}(i) : i \in \vSet, \, \rVal^{*}(i) = \bEnclose{\val^{*}(i)}^{-1}}$ must satisfy $\propM \rVal^{*} = \rVal^{*}$, where $\propM$ is the propagation matrix \eqref{eq:propMDef} for $\graph$.  But $\rVal^{*}$ must agree with $\rVal$, the reciprocal value vector for the limiting strategy, on the terminal nodes, since the values of those are pre-assigned.  It follows now from the proof of Lemma \ref{lem:limitingValuesTerm} that the non-terminal values $\rVal^{*}_{\ntl}$ must satisfy
\begin{equation*}
 \rVal^{*}_{\ntl} = \pEnclose{\id - \propMntnt}^{-1} \propMntt \rVal_{\tl},
\end{equation*}
and therefore $\rVal^{*} = \rVal$.  This completes the proof.
\end{proof}

For what terminating graphs $\graph$ is the Path Guessing Game fair, in the sense that the Guesser's expected fortune at the end of the game is equal to the \$1 that she started out with?

\begin{theorem} \label{thm:fairTerm}
Let $\graph$ be a terminating graph.  Then the Path Guessing Game on $\graph$ is fair if and only if each terminal node has a value of $1$ and each non-terminal node has out-degree at least $2$.
\end{theorem}
\begin{proof}
The game will be fair if and only if the value of every node is $1$.  From equation \eqref{eq:rValntl} this is equivalent to $\ones = \pEnclose{\id - \propMntnt}^{-1} \propMntt \ones$, and it follows from this that for each non-terminal $i$, the $\ith{i}$ row sum of $\propMntnt$ plus the $\ith{i}$ row sum of $\propMntt$ must equal $1$.  By the construction of the propagation matrix, this can happen if and only if $\graph$ has no non-terminal node of out-degree one.
\end{proof}

\section{Strongly Connected Graphs and the Discounted Game} \label{s:stronglyConnected}
%The rules for the Path Guessing Game on a strongly connected graph $\graph$ are identical to those for a terminating graph, the only difference being that there are no terminal nodes in a strongly connected graph.  The dynamics of the game are quite different, however.  Most obviously, the Chooser does not have the option of stopping the game by moving to a terminal node, so the game will have infinite duration.  Since the game lasts forever, there is also the possibility that the Guesser's fortune will grow without bound.

%In Section \ref{s:terminating}, optimal play on terminating graphs was derived by computing the values of the non-terminal nodes in terms of the pre-assigned values of the terminal nodes.  For strongly connected graphs, there are no terminal nodes, and hence no pre-assigned values.  We will nevertheless be able to assign values to the nodes based on limiting values for a truncated game, and then derive optimal play in terms of these limiting values.

Let $\graph = \pEnclose{\vSet,\eSet}$ be a strongly connected digraph on $\nNodes$ vertices.  The rules for the Path Guessing Game on $\graph$ are identical to those for a terminating graph, the only difference being that there are no terminal nodes in a strongly connected graph.  As in the previous section, we will derive optimal play on $\graph$ by considering the truncated game on $\graph$ obtained by the stopping the Path Guessing Game after $\tmNds$ steps.
%Similar to the truncated game on terminating graphs, the $\tmNds$-step game on a strongly connected graph is completely described by the $\nNodes$ path trees of length $\tmNds$ emanating from vertices.  The Guesser's fortune at the end of the game is simply her fortune after arriving at a leaf node of one of the path trees.  Thus, leaf nodes of the path trees are assigned the value of $1$.  Optimal play for the $\tmNds$-step game is now given by Theorem \ref{thm:optPlayTree}.

We are interested in the values of the root nodes of the path trees for the $\tmNds$-step game, as $\tmNds$ goes to infinity.  
%As before, we denote the value of the root node of the path tree of length $\tmNds$ emanating from node $i$ by $\val_{\tmNds}(i)$; the reciprocal of this value is denoted by $\rVal_{\tmNds}(i)$.  
Consider first the one-step game on $\graph$.  Equation \eqref{eq:treeValues} of Theorem \ref{thm:optPlayTree} holds for all vertices $i$, where the values of all the leaf nodes are $1$.  This can be summarized by the matrix equation
\begin{equation*}
 \rVal_{1} = \propM \ones,
\end{equation*}
where $\rVal_{1}$ is the vector of reciprocal values, $\ones$ is the vector of all ones, and $\propM = \pEnclose{\propE_{i,j} : i,j \in \vSet}$ is the $\nNodes \times \nNodes$ propagation matrix given by
\begin{equation} \label{eq:propMDefSC}
 \propE_{i,j} = 
 \begin{cases}
  1 / 2           \quad &\text{if $\outDeg_{i} = 1$ and $i \leadsTo j$,} \\
  1 / \outDeg_{i} \quad &\text{if $\outDeg_{i} \ge 2$ and $i \leadsTo j$,} \\
  0               \quad &\text{otherwise},
 \end{cases}
\end{equation}
where $\outDeg_{i}$ is the out-degree of vertex $i$.
Inductively, reciprocal values for the root nodes of the path trees for the $\tmNds$-step game are given by
\begin{equation*}
 \rVal_{\tmNds} = \propM^{\tmNds} \ones.
\end{equation*}

In contrast to the case of terminating graphs, limiting values for the $\tmNds$-step game may be infinite.  For instance, if $\graph$ is a $2$-cycle, then $\val_{\tmNds}(i) = 2^{\tmNds}$ for $i = 1, 2$.  However, recall from Theorem \ref{thm:optPlayFan} that optimal strategies depend only on the \emph{ratios} of values.  We will be able to show that limiting ratios continue exist in the strongly connected case, and that the corresponding strategies are optimal.

Towards this end, consider the following \emph{discounted game} on the strongly connected graph $\graph$:  play proceeds just as before, but after each payoff, the Guesser's fortune is multiplied by a fixed \emph{discount factor} $\dFactor \in \ocint{0}{1}$.  Intuitively, one can think of the discounted game as modeling a situation in which the real value of money is decreasing with time, as in an inflationary economy.  (In this case, the reciprocal of the discount rate would be the inflation rate.)  It is clear that optimal strategies for the discounted game are the same as for the undiscounted game.  

For the discounted $1$-step game, the values of all the leaf nodes are $\dFactor$.  Reciprocal values are therefore $1/\dFactor$, and so the reciprocal values of the root nodes in the discounted game are
\begin{equation*}
 \rValD_{1} = \propM \pEnclose{ \dFactor^{-1} \ones} = \dFactor^{-1} \propM \ones.
\end{equation*}
Inductively, values for the root nodes of the path trees for the discounted $\tmNds$-step game are
\begin{equation} \label{eq:valPropSC}
 \rValD_{\tmNds} = \dFactor^{-\tmNds} \propM \ones,
\end{equation}
and the idea is to determine the value of the discount factor $\dFactor$ that makes the limiting reciprocal values finite and nonzero.  

It turns out that the correct choice is to make the discount factor equal to the largest positive eigenvalue of the propagation matrix $\propM$.  Before proving this, we will make one more assumption on the graph $\graph$.  In addition to being strongly connected, we will assume that $\graph$ is \emph{aperiodic}.  For a strongly connected graph, aperiodicity means that the greatest common divisor of the lengths of all the cycles in $\graph$ is one.  (Note that any strongly connected graph that contains a loop is aperiodic.)  Aperiodicity rules out certain cyclic phenomena that, while not unduly hard to characterize, serve mainly to cloud the important issues.  
The next lemma collects the properties of the propagation matrix $\propM$ that we will need.
\begin{lemma} \label{lem:propMPropertiesSC}
Let $\graph$ be a strongly connected aperiodic digraph, with associated propagation matrix $\propM$ given by \eqref{eq:propMDefSC}.  Then we have the following:
\begin{enumerate}
 \item \label{it:propMProp1} $\propM$ has a positive maximal eigenvalue $\propMEV$, with the property that any other eigenvalue $\lambda$ of $\propM$ satisfies $\abs{\lambda} < \propMEV$.
 \item \label{it:propMProp2} There are strictly positive right and left eigenvectors $\rEV$ and $\lEV$ respectively associated with the maximal eigenvalue $\propMEV$.
 \item \label{it:propMProp3} The right and left eigenspaces of $\propM$ associated with $\propMEV$ (and containing $\rEV$ and $\lEV$ respectively) each have dimension $1$.
 \item \label{it:propMProp4} No other eigenvector of $\propM$ is positive.
 \item \label{it:propMProp5} The maximal eigenvalue $\propMEV$ satisfies $\frac{1}{2} \le \propMEV \le 1$.
\end{enumerate}
\end{lemma}
\begin{proof}
Since $\graph$ is strongly connected and aperiodic, it follows that $\propM$ is irreducible and primitive, and properties \ref{it:propMProp1} through \ref{it:propMProp4} follow from the standard Perron-Frobenius theory of nonnegative matrices (see \cite{Minc}, for instance).  Property \ref{it:propMProp5} follows from the fact that all the row sums of $\propM$ are between one half and one.
\end{proof}

We are now in a position to prove that limiting values exist for the discounted game.
\begin{lemma}[Limiting Values for the Discounted Game] \label{lem:limitingValuesSC}
Let $\graph$ be a strongly connected aperiodic digraph, with associated propagation matrix $\propM$ given by \eqref{eq:propMDefSC}.  Let $\propMEV$ be the maximal eigenvalue of $\propM$, with associated positive right and left eigenvectors $\rEV$ and $\lEV$ respectively.  Then we have
\begin{equation} \label{eq:propMLimitSC}
 \lim_{\tmNds \approaches \infty} \propMEV^{-\tmNds} \propM^{\tmNds} = \frac{\rEV \lEV^{\transpose}}{\rEV^{\transpose} \lEV},
\end{equation}
where the limit is a strictly positive matrix.  
Moreover, the limiting reciprocal value vector $\rValD$ for the discounted game (with discount factor $\dFactor = \propMEV$) given by
\begin{equation} \label{eq:limitingValuesSC}
 \rValD
 = \lim_{\tmNds \approaches \infty} \propMEV^{-\tmNds} \propM^{\tmNds} \ones 
 = \frac{\rEV \lEV^{\transpose}}{\rEV^{\transpose} \lEV} \ones
\end{equation}
is a positive right eigenvector of $\propM$ corresponding to $\propMEV$:
\begin{equation}
 \propM \rValD = \propMEV \rValD
\end{equation}
\end{lemma}
\begin{proof}
Equation \eqref{eq:propMLimitSC} is easily derived from Lemma \ref{lem:propMPropertiesSC} and the Jordan canonical form for $\propM$ (see \cite{Lancaster}, for instance).  Equation \eqref{eq:limitingValuesSC} follows directly from \eqref{eq:valPropSC} and \eqref{eq:propMLimitSC}.  Finally,
\begin{align*}
 \propM \rValD
  &= \propM \lim_{\tmNds \approaches \infty} \propMEV^{-\tmNds} \propM^{\tmNds} \ones \\
  &= \propMEV \lim_{\tmNds \approaches \infty} \propMEV^{-(\tmNds+1)} \propM^{\tmNds+1} \ones \\
  &= \propMEV \rValD,
\end{align*}
and $\rValD$ is positive because both $\rEV$ and $\lEV$ are.
\end{proof}

Lemma \ref{lem:limitingValuesSC} shows that setting the discount factor $\dFactor$ equal to the maximal eigenvalue $\propMEV$ of the propagation matrix results in a discounted game in which the limiting values are all finite and nonzero.  
%Moreover $\dFactor = \propMEV$ is the unique discount factor that does this: if $\dFactor > \propMEV$, then the limiting values are all infinite, while if $\dFactor < \propMEV$ the limiting values are all $0$.  
From this point onwards, when we refer to the ``discounted game'', it is understood that the discount factor is the maximal eigenvalue of $\propM$.  Also, for ease of notation we drop the tilde and simply use $\rVal$ to refer to the limiting reciprocal values for the discounted game.

Turning towards the non-truncated infinite duration game, we now explore strategies corresponding to the limiting values.   For the Chooser, the transition probabilities are
\begin{equation} \label{eq:cOptProbSC0}
 \cProb_{i,j} = \cPr{\pos_{\tmNdx+1} = j}{\pos_{\tmNdx} = i} = 
 \begin{cases}
  1 &\quad \text{if $\outDeg_i = 1$} \\
  \rVal_{j} / \pEnclose{\sum_{k : i \leadsTo k} \rVal_{k}} &\quad \text{if $\outDeg_i \ge 2$}
 \end{cases}
\end{equation}
while the Guesser has guessing probabilities
\begin{equation} \label{eq:gOptProbSC}
 \gProb_{i,j} = \cPr{\gGuess_{\tmNdx+1} = j}{\pos_{\tmNdx} = i} = 
 \begin{cases}
  1 &\quad \text{if $\outDeg_i = 1$} \\
  \pEnclose{\cProb_{i,j} - \gParam \cProb_{i,\min}} / \pEnclose{1 - \outDeg_{i} \gParam \cProb_{i,\min}} &\quad \text{if $\outDeg_i \ge 2$}
 \end{cases}
\end{equation}
and wagers
\begin{equation} \label{eq:gOptWagerSC}
 \wager_{i} = 
 \begin{cases}
  1 & \quad \text{if $\outDeg_{i} = 1$} \\
  1 - \outDeg_{i} \gParam \cProb_{i,\min} & \quad \text{if $\outDeg_{i} \ge 2$} 
 \end{cases}
% \cPr{\gWager_{\tmNdx+1} = 1 - \outDeg_{i} \gParam \cProb_{i,\min}}{\pos_{\tmNdx} = i} = 1 
\end{equation}
when at node $i$.  Here again $\gParam$ is arbitrary in $\ocint{0}{1}$, and $\cProb_{i,\min}$ is the probability of the vertex least likely to be chosen by the Chooser:
\begin{equation*}
 \cProb_{i,\min} = \min \cEnclose{\cProb_{i,j} : \text{$j \in \vSet$, $i \leadsTo j$}}
\end{equation*}
We will refer to the strategies given by \eqref{eq:cOptProbSC0}, \eqref{eq:gOptProbSC}, and \eqref{eq:gOptWagerSC} as the \emph{limiting strategies} for the players.  The next theorem gives the basic properties of the game under the limiting strategies.
\begin{theorem}[Limiting Strategies] \label{thm:limitingStrategiesSC}
Let $\graph = \pEnclose{\vSet,\eSet}$ be a strongly connected aperiodic graph.  Let $\propM$ be the propagation matrix for $\graph$ defined by \eqref{eq:propMDefSC}.  Let $\propMEV$ be the maximal eigenvalue of $\propM$, and let $\val = \pEnclose{\val_{i} : i \in \vSet}$ be the vector of limiting values.  Then under the limiting strategies \eqref{eq:cOptProbSC0}, \eqref{eq:gOptProbSC}, and \eqref{eq:gOptWagerSC} the Chooser's transition probability matrix $\cProbM$ is diagonally similar to the propagation matrix:
\begin{equation} \label{eq:cOptProbMSC}
 \cProbM = \frac{1}{\propMEV} \valM \propM \valM^{-1},
\end{equation}
where $\valM = \text{diag}(\val)$ is the diagonal matrix of limiting values.  Under the limiting strategies, the Guesser's fortune satisfies
\begin{equation} \label{eq:eFortSC}
 \lim_{\tmNdx \approaches \infty} \propMEV^{\tmNdx} \cE{\fort_{\tmNdx}}{\pos_{0} = i} = c \, \val_{i},
\end{equation}
where $c$ is a positive constant that depends only on $\graph$.
\end{theorem}
\begin{proof}
Because the limiting reciprocal value vector $\rVal$ satisfies $\propM \rVal = \propMEV \rVal$, equation \eqref{eq:cOptProbSC0} reduces to
\begin{equation*} 
\cProb_{i,j} = \frac{1}{\propMEV} \, \propE_{i,j} \, \rVal_{j} / \rVal_{i} = \frac{1}{\propMEV} \, \propE_{i,j} \, \val_{i} / \val_{j}
\end{equation*}
for all $i$ and $j$, and equation \eqref{eq:cOptProbMSC} follows.

For \eqref{eq:eFortSC}, start by considering a vertex $i$ whose out-degree $\outDeg_{i}$ is at least $2$.  Then for any adjacent vertex $j$ we have
\begin{align*}
 \cE{\fort_{1}}{\pos_{0} = i, \pos_{1} = j} 
  &= \pEnclose{1 + (\outDeg_{i} - 1) \wager_{i}} \gProb_{i,j} + \pEnclose{1 - \wager_{i}} \pEnclose{1 - \gProb_{i,j}} \\
  &= \wager_{i} \pEnclose{\outDeg_{i} \gProb_{i,j} - 1} + 1 \\
  &= \outDeg_{i} \pEnclose{\cProb_{i,j} - \gParam \cProb_{i,\min}} - \pEnclose{1 - \outDeg_{i} \gParam \cProb_{i,\min}} + 1 \\
  &= \outDeg_{i} \cProb_{i,j} \\
  &= \frac{1}{\propMEV} \frac{\val_{i}}{\val_{j}}.
\end{align*}
Therefore
\begin{equation} \label{eq:valPropOneStepSC}
 \val_{i} = \cE{\fort_{1}}{\pos_{0} = i, \pos_{1} = j} \propMEV \, \val_{j}
\end{equation}
In fact, \eqref{eq:valPropOneStepSC} holds for vertices $i$ of out-degree $1$ as well.  To see this, note that if $i$ has out-degree $1$, and $i \leadsTo j$, then by equation \eqref{eq:propMDefSC} and $\propM \rVal = \propMEV \rVal$ we see that 
\begin{equation*} %\label{eq:mark1}
\val_{i} = 2 \propMEV \val_{j}.
\end{equation*}
But since $i$ has out-degree $1$, $\fort_{1} = 2$ with probability $1$, so in fact the last equation is equivalent to \eqref{eq:valPropOneStepSC}.  Thus, \eqref{eq:valPropOneStepSC} holds for \emph{all} vertices $i$ and $j$ such that $i \leadsTo j$.
%Next, note that if $i \leadsTo j \leadsTo k$, then by the Markov property
%\begin{align*}
% \cE{\fort_{2}}{\pos_{0} = i, \pos_{1} = j, \pos_{2} = k, \fort_{1}} 
%  &= \fort_{1}  \cE{\fort^{\prime}_{1}}{\pos_{0} = j, \pos_{1} = k} \\
%  &= \fort_{1} \frac{1}{\propMEV} \frac{\val_{j}}{\val_{k}}.
%\end{align*}
%Thus
%\begin{align*}
% \cE{\fort_{2}}{\pos_{0} = i, \pos_{1} = j, \pos_{2} = k} 
%  &= \frac{1}{\propMEV} \frac{\val_{j}}{\val_{k}} \cE{\fort_{1}}{\pos_{0} = i, \pos_{1} = j} \\
%  &= \frac{1}{\propMEV} \frac{\val_{j}}{\val_{k}} \, \frac{1}{\propMEV} \frac{\val_{i}}{\val_{j}} \\
%  &= \frac{1}{\propMEV^{2}} \frac{\val_{i}}{\val_{k}}
%\end{align*}
%Since this is true for all $j$ such that $i \leadsTo j$, we have
%\begin{equation*}
% \cE{\fort_{2}}{\pos_{0} = i, \pos_{2} = k} = \frac{1}{\propMEV^{2}} \frac{\val_{i}}{\val_{k}}
%\end{equation*}
%whenever $k$ is reachable from $i$ in exactly $2$ steps; equivalently
%\begin{equation*} \label{eq:valPropTwoStepSC}
% \val_{i} = \cE{\fort_{2}}{\pos_{0} = i, \pos_{1} = j} \propMEV^{2} \, \val_{j}.
%\end{equation*}
%Inductively we see that 
%\begin{equation} \label{eq:valPropTStepSC}
% \val_{i} = \cE{\fort_{\tmNdx}}{\pos_{0} = i, \pos_{\tmNdx} = j} \propMEV^{\tmNdx} \, \val_{j}
%\end{equation}
%holds for all $\tmNdx$ and all $j$ that are reachable from $i$ in exactly $\tmNdx$ steps.  
Now using the Markov property we have
\begin{align*}
 \cE{\fort_{\tmNdx}}{\pos_{0} = i} 
  &= \sum_{j} \cE{\fort_{\tmNdx}}{\pos_{0} = i, \pos_{1} = j} \cProb_{i,j} \\
  &= \sum_{j} \cE{\fort_{1}}{\pos_{0} = i, \pos_{1} = j} \cE{\fort^{\prime}_{\tmNdx-1}}{\pos_{0} = j} \cProb_{i,j} \\
  &= \sum_{j} \frac{1}{\propMEV} \frac{\val_{i}}{\val_{j}} \cE{\fort_{\tmNdx-1}}{\pos_{0} = j} \cProb_{i,j} \\
  &= \sum_{j} \frac{1}{\propMEV^{2}} \frac{\val_{i}^{2}}{\val_{j}^{2}} \propE_{i,j} \cE{\fort_{\tmNdx-1}}{\pos_{0} = j}
\end{align*}
The matrix form of this equation is
\begin{equation*}
 \zeta_{\tmNdx} = \propMEV^{-2} \valM^{2} \propM \valM^{-2} \zeta_{\tmNdx-1}.
\end{equation*}
where $\zeta_{\tmNdx}(i) = \cE{\fort_{\tmNdx}}{\pos_{0} = i}$.  Since $\zeta_{0} = \ones$, we have
\begin{equation} \label{eq:eFortSC2}
 \zeta_{\tmNdx} = \propMEV^{-2 \tmNdx} \valM^{2} \propM^{\tmNdx} \valM^{-2} \ones,
\end{equation}
and so by Lemma \ref{lem:limitingValuesSC}
\begin{align*}
 \lim_{\tmNdx \approaches \infty} \propMEV^{\tmNdx} \zeta_{\tmNdx} 
  &= \propMEV^{-\tmNdx} \valM^{2} \propM^{\tmNdx} \valM^{-2} \ones \\
  &= \valM^{2} \, \, \frac{\rEV \lEV^{\transpose}}{\rEV^{\transpose} \lEV} \, \, \valM^{-2} \ones.
\end{align*}

To see that the limit is a multiple of the limiting value vector $\val$, denote $G = \propMEV^{-2} \valM^{2} \propM \valM^{-2}$, so that
\begin{equation} \label{eq:Geq}
 \zeta_{\tmNdx} = G^{\tmNdx} \ones.
\end{equation}
Note that $G$ has the same pattern of zero and non-zero entries as $M$, and hence is irreducible and primitive.  We claim that $1/\propMEV$ is the maximal eigenvalue of $G$, and that the limiting value vector $\val$ is an associated eigenvector.  To see this, first observe that 
\begin{align*}
 G \val
  &= \propMEV^{-2} \valM^{2} \propM \valM^{-2} \val \\
%  &= \propMEV^{-2} \valM^{2} \propM \valM^{-1} \ones \\
  &= \propMEV^{-2} \valM^{2} \propM \rVal \\
  &= \propMEV^{-2} \valM^{2} r \rVal \\
%  &= \propMEV^{-1} \valM \ones \\
  &= \propMEV^{-1} \val,
\end{align*}
so that $1/\propMEV$ is indeed an eigenvalue of $G$, with associated eigenvector $\val$.  To see that $1/\propMEV$ is the maximal eigenvalue of $G$, note that if $\pEnclose{\lambda, x}$ is any eigenpair for $G$, then  
$G x = \lambda x$ implies that $\propM (\valM^{-2} x) = \propMEV^{2} \lambda (\valM^{-2} x)$, so that $\propMEV^{2} \lambda$ is an eigenvalue of $\propM$, with associated eigenvector $\valM^{-2} x$.  But $\propMEV$ is the \emph{maximal} eigenvalue of $\propM$, so we must have $\abs{\propMEV^{2} \lambda} \le \propMEV$, which implies $\abs{\lambda} \le 1/\propMEV$, which means that $1/\propMEV$ is the maximal eigenvalue of $G$.  Now since $G$ is primitive with maximal eigenvalue $1/\propMEV$, it follows that $\propMEV^{\tmNdx} G^{\tmNdx} \ones$ converges to some multiple of the associated eigenvector $v$.  By \eqref{eq:Geq} this means that
\begin{equation*}
 \lim_{\tmNdx \approaches \infty} \propMEV^{\tmNdx} \cE{\fort_{\tmNdx}}{\pos_{0} = i} = c \, \val_{i}.
\end{equation*}
This completes the proof.
\end{proof}

Because $\graph$ is strongly connected and aperiodic, the random walk on $\graph$ that arises from the limiting strategies is ergodic, and has a unique invariant measure $\invMsr$.  Under this invariant measure, the \emph{discounted} fortune process
\begin{equation} %\label{eq:discFortProcess}
 \dFort_{\tmNdx} = \propMEV^{\tmNdx} \fort_{\tmNdx}
\end{equation}
is stationary.  Let $\ssE{\cdot}$ denote the expectation operator under the invariant measure (i.e. assuming that the initial position $\pos_{0}$ of the players has distribution $\invMsr$).  The \emph{steady-state fortunes} defined by
\begin{equation*}
 \eta_{j} = \sscE{\dFort_{\tmNdx}}{\pos_{\tmNdx} = j}, \quad j \in \vSet
\end{equation*}
represent the Guesser's average discounted fortune when located at vertex $j$.  By stationarity, they do not depend on $\tmNdx$.  The next theorem characterizes the invariant measure and steady-state fortunes.

\begin{theorem} \label{thm:markovPropsSC}
Let $\graph$ be a strongly connected aperiodic graph with propagation matrix $\propM$ and associated maximal eigenvalue $\propMEV$.  Then the invariant measure $\invMsr$ for the position process $\pEnclose{\pos_{\tmNdx} : \tmNdx \ge 0}$ from the limiting strategies is the entry-wise product of appropriately scaled right and left eigenvectors $\rEV$ and $\lEV$ respectively of $\propM$ associated with $\propMEV$:
\begin{equation} \label{eq:invMsrSC}
 \invMsr = \frac{1}{\rEV^{\transpose} \lEV} \pEnclose{\rEV_{i} \lEV_{i} : i \in \vSet}.
\end{equation}
The steady state fortunes are of the form
\begin{equation} \label{eq:ssForts}
  \sscE{\dFort_{\tmNdx}}{\pos_{\tmNdx} = j} = c \, \frac{\invMsr_{j}}{\val_{j}},
\end{equation}
where $c$ is a constant that depends on $\graph$.
\end{theorem}
\begin{proof}
The invariant measure $\invMsr$ satisfies $\cProbM^{\transpose} \invMsr = \invMsr$.  But $\cProbM = \propMEV^{-1} \valM \propM \valM^{-1}$, so the invariant measure must satisfy $\propM^{\transpose} \valM \invMsr = \propMEV \valM \invMsr$.  By primitivity, $\valM \invMsr = c \lEV$, where $\lEV$ is the (unique) left eigenvector of $\propM$ correponding to $\propMEV$.  Thus 
\begin{equation*}
\invMsr = c \valM^{-1} \lEV = c \pEnclose{\rVal_{i} \lEV_{i} : i \in \vSet},
\end{equation*}
and \eqref{eq:invMsrSC} follows because the reciprocal value vector $\rVal$ is a right eigenvector of $\propM$.

For \eqref{eq:ssForts}, by conditioning on the position of the players at time $\tmNdx - 1$, one sees that
\begin{align*}
 \sscE{\dFort_{\tmNdx}}{\pos_{\tmNdx} = j} 
  &= \sum_{i : i \leadsTo j} \frac{\val_{i}}{\val_{j}} \sscE{\dFort_{\tmNdx}}{\pos_{\tmNdx} = i} \cProb_{i,j} \\
  &= \sum_{i : i \leadsTo j} \frac{1}{\propMEV} \frac{\val_{i}^{2}}{\val_{j}^{2}} \propE_{i,j} \sscE{\dFort_{\tmNdx}}{\pos_{\tmNdx} = i}.
\end{align*}
The matrix form of this equation is
\begin{equation*}
 \propMEV \eta^{\transpose} = \eta^{\transpose} \valM^{2} \propM \valM^{-2},
\end{equation*}
where $\eta_{j} = \sscE{\dFort_{\tmNdx}}{\pos_{\tmNdx} = j}$.  By the primitivity of $\propM$ we have $\valM^{2} \eta = c \lEV$, and \eqref{eq:ssForts} now follows.
\end{proof}

\begin{theorem}[Optimal play on strongly connected graphs] \label{thm:optPlaySC}
Let $\graph$ be a strongly connected aperiodic graph.  Then the limiting strategies \eqref{eq:cOptProbMSC}, \eqref{eq:gOptProbSC}, and \eqref{eq:gOptWagerSC} are optimal for the infinite-duration Path Guessing Game on $\graph$.
\end{theorem}
\begin{proof}[Proof Sketch]
As in the proof of Theorem \ref{thm:optPlayTerm}, we can restrict attention to purely positional strategies.  For any choice of purely positional strategies by the Players, there will exist some discount factor $\dFactor$ such that the limits
\begin{equation*}
 \lim_{\tmNds \approaches \infty} \dFactor^{\tmNdx} \cE{\fort_{\tmNdx}}{\pos_{0} = i} = \val^{\star}_{i}
\end{equation*}
exist and are non-zero.  For this choice of strategies, the discounted game starting at node $i$ is therefore equivalent to the game on the fan with root node $i$, whose leaves are the set of vetices $j$ such that $i \leadsTo j$ in $\graph$, and with vertex $j$ having the value $\dFactor \, \val^{\star}_{j}$.   By the results of Section \ref{s:fans}, the players are playing optimally if and only if $\propM \rVal^{\star} = \dFactor \rVal^{\star}$, where $\rVal^{\star}$ is the vector of reciprocals of the entries of $\val^{\star}$.  But $\propM$ is irreducible and primitive, so its \emph{only} positive eigenvector is the reciprocal value vector $\rVal$ of Lemma \ref{lem:limitingValuesSC}.  Thus $\rVal^{\star} = \rVal$, and it follows that the optimal strategies are the same as the limiting strategies.
\end{proof}

For what strongly connected aperiodic graphs $\graph$ is the Path Guessing Game fair, in the sense that the Guesser's expected fortune at the end of the game is equal to the \$1 that she started out with?
\begin{theorem} \label{thm:fairSC}
Let $\graph$ be a strongly connected aperiodic graph.  Then the Path Guessing Game on $\graph$ is fair if and only if each vertex of $\graph$ has out-degree at least two.
\end{theorem}
\begin{proof}
If each vertex has out-degree two or more, then the propagation matrix $\propM$ is stochastic, implying that $\propMEV = 1$ and $\rVal = \ones$, and it follows from \eqref{eq:eFortSC2} that the game is fair: $\cE{\fort_{\tmNdx}}{\pos_{0} = i} = 1$ for all $i \in \vSet$.  On the other hand, if there exists a vertex $i$ of out-degree one, then $\propM$ is strictly sub-stochastic, and in fact the $\ith{i}$ component of $\propM \ones$ equals $1/2$.  By primitivity there exists an integer $\tmNdx_{0} > 0$ such that $\propM^{\tmNdx_{0}}$ is strictly positive.  It follows that that \emph{all} the components of $\propM^{\tmNdx_{0} + 1} \ones$ are strictly less than one, and hence $\lim_{\tmNdx \approaches \infty} \propM^{\tmNdx} \ones = 0$.  Therefore the maximal eigenvalue $\propMEV$ of $\propM$ must be strictly less than one, and therefore the Guesser's expected fortune under optimal play grows without bound.
\end{proof}

\section{The Lying Oracle Game} \label{s:lyingOracle}
In this section we apply the previous results to the \emph{Lying Oracle Game}.  As described in \cite{Koether}, this is a two-player game between an Oracle and a Bettor.  The game is based on a sequence of $n$ tosses of a fair coin.  The Oracle has perfect foreknowledge of each toss, but may choose to lie to the Bettor in up to $k$ of the $n$ tosses.  Before each toss, the Bettor privately writes down a guess as to whether she thinks the Oracle will lie or tell the truth on this toss.  The Bettor also announces a wager on her guess.  The Oracle then announces his (possibly mendacious) prediction.  Finally the coin is then tossed, and the outcome is recorded.  If the Bettor was correct, she wins the amount of her wager; if she was incorrect, she loses the amount of her wager.  Play continues in this fashion until the $n$ of tosses are completed.

The Lying Oracle Game has been shown to have a reciprocal relationship with the continuous version of \emph{Ulam's search problem} \cite{Ravikumar}, which has been used as a model of binary search in the presence of errors \cite{Rivest}.  In Ulam's search problem, a questioner is searching for a number in the interval $\ccint{0}{1}$ chosen by a responder.  The questioner asks the responder $n$ questions about the number's location in the interval, but the responder may lie up to $k$ times.  The questioner seeks to find a subset of smallest measure that contains the chosen number, while the responder seeks to maximize the measure of that subset.  Under optimal play, the measure of the questioner's subset is the reciprocal of the Bettor's fortune in the Lying Oracle Game.

The Lying Oracle Game can be generalized in several ways.  In \cite{Koether2} the authors consider the game when the coin is not fair.   The ``at most $k$ lies in $n$ tosses'' rule can be generalized to a set of arbitrary ``lying patterns''.  Here, we make the point that the original Lying Oracle Game (with a fair coin) can be seen as a special case of the Path Guessing Game.  Moreover, the results from the previous sections can be used to derive optimal play in the Lying Oracle Game when the number of coin tosses is infinite.

\begin{example}
Consider the game $\graph_{n,1}$ in which the Oracle can lie at most one time in any block of $n$ predictions.  
The corresponding graph is shown in Figure 1.

\begin{center}
 \includegraphics[width=2in]{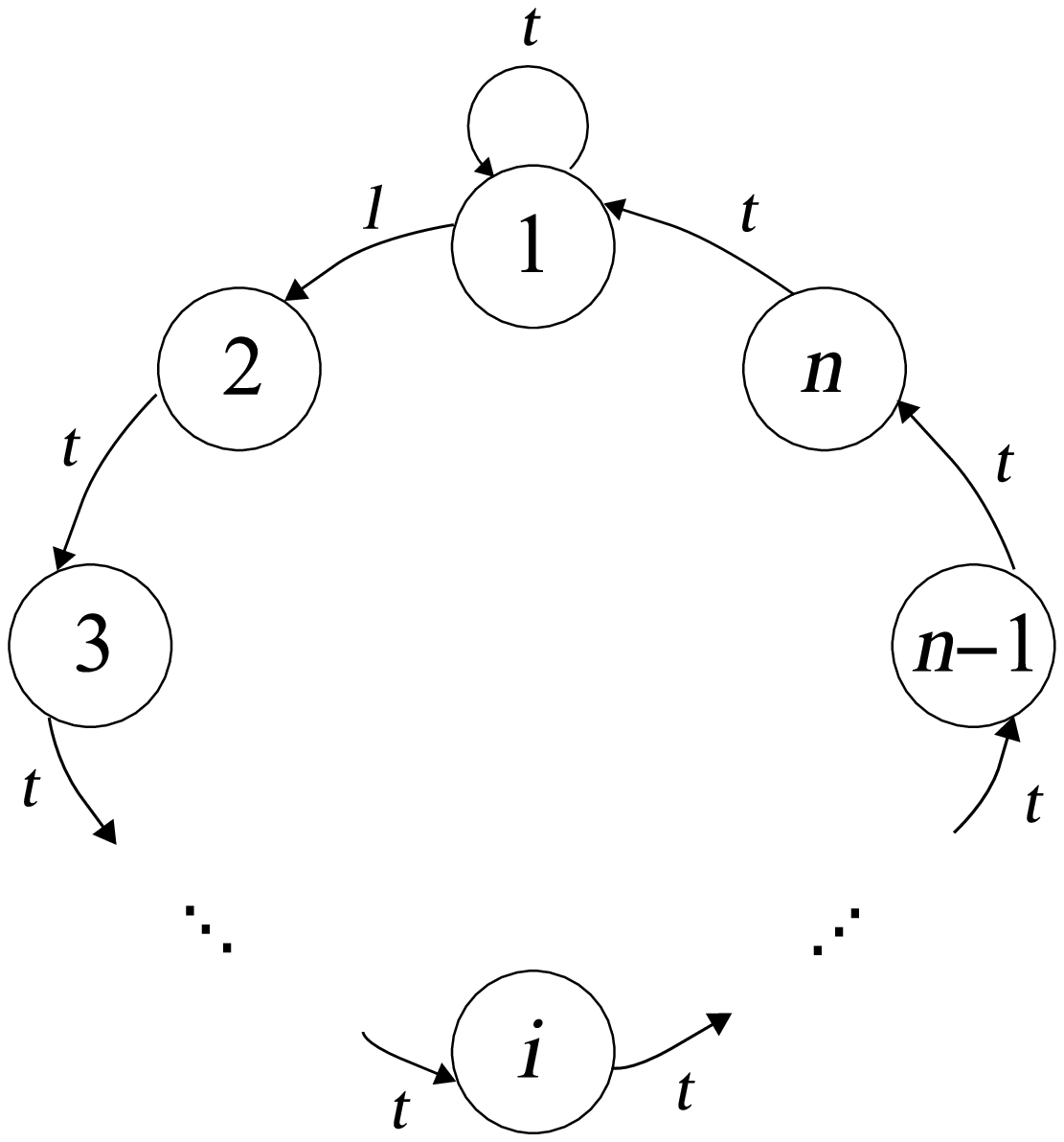}
\end{center}
%\begin{center}
Figure 1.  The game $\graph_{n,1}$.  The Oracle can lie at most $1$ time in any block of $n$ tosses.  Play starts at vertex $1$.
%\end{center}
\vskip 0.5cm

Note that since the out-degree of each vertex is two or less, our payoff rule \eqref{eq:payoffRule} coincides with the payoff rule for the Lying Oracle Game.  The propagation matrix is $\propM = \frac{1}{2} \adjM$, where 
\begin{equation*}
 \adjM = 
 \begin{pmatrix}
  1 & 1 & 0 & 0 & \dots & 0 \\
  0 & 0 & 1 & 0 & \dots & 0 \\
  0 & 0 & 0 & 1 & \dots & 0 \\
  \vdots & \vdots & \vdots & \vdots & \ddots & \vdots \\
  0 & 0 & 0 & 0 & \dots & 1 \\
  1 & 0 & 0 & 0 & \dots & 0
 \end{pmatrix}
\end{equation*}
The maximal eigenvalue of $\propM$ is $\propMEV = \frac{1}{2} \adjMEV$, where $\adjMEV = \adjMEV(n)$ is the largest positive solution of
\begin{equation*} \label{eq:one-in-n-charPoly}
 0 = \det \pEnclose{\adjMEV \id - \adjM} = \adjMEV^{n} - \adjMEV^{n-1} - 1.
\end{equation*}
The row-sums of $\adjM$ tell us that $\adjMEV \in \ccint{1}{2}$, and it is not hard to show that as $n$ approaches infinity $\adjMEV$ decreases to $1$ monotonically.  Thus $\propMEV$ approaches $1/2$, as we would expect.  The right eigenvector of $\propM$ corresponding to $\propMEV$ is
\begin{equation*}
 \rEV = \pEnclose{\adjMEV^{n-1}, 1, \adjMEV, \adjMEV^{2}, \dots, \adjMEV^{n-3}, \adjMEV^{n-2}},
\end{equation*}
and the left eigenvector corresponding to $\propMEV$ is
\begin{equation*}
 \lEV = \pEnclose{\adjMEV^{n-1}, \adjMEV^{n-2}, \adjMEV^{n-3}, \dots, \adjMEV, 1}.
\end{equation*}
The Oracle's optimal strategy at node $1$ (the only non-trivial node) is
\begin{equation*}
 \cProb_{1,1} = \adjMEV^{-1}, \quad \cProb_{1,2} = \adjMEV^{-n}.
\end{equation*}
In other words, the Oracle tells the truth at node $1$ with probability $\adjMEV^{-1}$, and lies with probability $1 - \adjMEV^{-1} = \adjMEV^{-n}$.
Note that as $n$ increases without bound, the probability that the Oracle tells the truth when he is at node $1$ approaches $1$:
\begin{equation*}
 \lim_{n \approaches \infty} \cProb_{1,1} = \lim_{n \approaches \infty} \adjMEV^{-1} = 1
\end{equation*}
The Bettor's minimum risk ($\gParam = 1$) optimal strategy is
\begin{align*}
 &\gProb_{1,1} = 1, \quad \gProb_{1,2} = 0 \\
 &\wager_{1,1} = \wager_{1,2} = \adjMEV^{-1} - \adjMEV^{-n}
\end{align*}
Under optimal play, the players perform a random walk through the graph, with invariant measure
\begin{equation*}
 \invMsr = \frac{1}{\adjMEV^{n}+n-1} \pEnclose{\adjMEV^{n}, 1, 1, \dots, 1}
\end{equation*}
Observe that the fraction of time that the Oracle lies under optimal play is
\begin{equation*}
 \invMsr_{2} = \frac{1}{\adjMEV^{n}+n-1},
\end{equation*}
which increases to $1/n$ as $n$ approaches infinity.  Thus under optimal play the Oracle lies slightly less often than the rules allow.
%It is not hard to show that the maximum root $\adjMEV = \adjMEV_{n}$ of \eqref{e:one-in-n-charPoly} satisfies 
%\begin{equation*}
% \lim_{n \approaches \infty} \adjMEV_{n}^{n} = \infty.
%\end{equation*}
%Thus for large $n$, the players spend far more time at vertex $1$ than at any other single vertex.
\hfill \expSymb
\end{example}

We remark that if the Oracle's set of allowed ``lying patterns'' can be expressed in the form a finite set of finite \emph{forbidden patterns} - sequences of truths and lies that the Oracle must avoid -  then the Lying Oracle Game (with a fair coin) is equivalent to the Path Guessing Game on a certain finite graph.  All games in which the Oracle can lie at most $k$ times in any block of $n$ statements fall into this class, and as such can be analyzed using the results of Section \ref{s:stronglyConnected}.  Our approach also leads to new variants of the Lying Oracle Game in which the Oracle is allowed to end the game under \emph{certain conditions}, rather than after a specified number of tosses.   For instance, we can consider a game in which the Oracle can lie at most $k$ times in any block of $n$ tosses, and can stop the game after any toss on which he told the truth.  In these games, the Guesser tries to predict whether the Oracle will lie, tell the truth, or stop the game.  Such games correspond to a Path Guessing Game on a terminating graph.  

\begin{example}
Consider the game in which the Oracle can lie at most $1$ time in any block of $n$ tosses, and can stop the game on the first round, or after any round on which he told the truth.  The graph for this game is shown in Figure 2.

\begin{center}
 \includegraphics[width=2in]{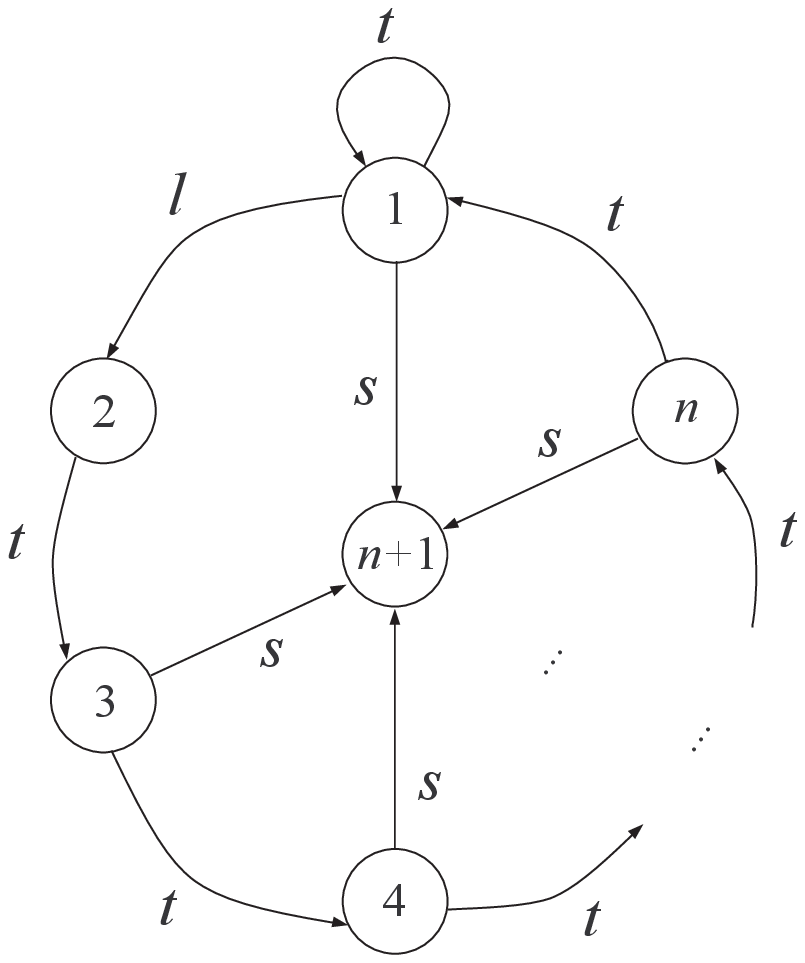}
\end{center}
%\begin{center}
Figure 2.  The Oracle can lie at most $1$ time in any block of $n$ tosses, and the Oracle can choose to stop after any toss on which he told the truth.  Play starts at vertex $1$.
%\end{center}
\vskip 0.5cm

There are $n$ nonterminal nodes and one terminal node.  Using the results of Section \ref{s:terminating}, we find that under optimal play the probability $\cProb_{i,n+1}$ that the Oracle stops the game from node $i$ is given by
\begin{equation*}
 \cProb_{i,n+1} = 
 \begin{cases}
  \frac {{2}^{n}-1}{3\,\pEnclose{{2}^{n-1}+{2}^{n-2}-1}} \quad &\text{if $i = 1$} \\
  0 \quad &\text{if $i = 2$} \\
  \frac {{2}^{n}-1}{{2}^{n+1}-{2}^{i-2}-2} \quad &\text{if $3 \le i \le n$}
 \end{cases}
\end{equation*}
It is easy to show that if $n$ and $i$ go to infinity concurrently such that $i/n$ approaches $x \in \ooint{0}{1}$, then $\cProb_{i,n+1}$ approaches $1/2$.  We also have $\cProb_{1,n+1}$ approaching $4/9$ as $n$ approaches infinity.  Thus, even for large $n$, the game duration is likely to be very short.
\hfill \expSymb
\end{example}

\paragraph*{Acknowledgement.}  The author thanks Robb Koether for many stimulating conversations on the Lying Oracle Game and its generalizations.

%\newpage

%\begin{center}
% \includegraphics[width=2in]{graph-one-in-n}
%\end{center}
%%\begin{center}
%Figure 1.  The game $\graph_{n,1}$.  The Oracle can lie at most $1$ time in any block of $n$ tosses.  Play starts at vertex $1$.
%%\end{center}
%

%\newpage

%\begin{center}
% \includegraphics[width=2in]{graph-1-in-n-term}
%\end{center}
%%\begin{center}
%Figure 2.  The Oracle can lie at most $1$ time in any block of $n$ tosses, and the Oracle can choose to stop after any toss on which he told the truth.  Play starts at vertex $1$.
%%\end{center}


\newpage

\begin{thebibliography}{5}

\bibitem{Koether}
Robb Koether and John Osoinach, Outwitting the Lying Oracle, \textsl{Mathematics Magazine}, 78 (2005), 98-109.

\bibitem{Koether2}
Robb Koether, Marcus Pendergrass,  and John Osoinach, The Lying Oracle with a Biased Coin, to appear in \textsl{Journal of Applied Probability}.

\bibitem{Ravikumar}
B. Ravikumar, Some connections between the lying oracle problem and Ulam's search problem, \textsl{Proceedings of AWOCA 2005, the Sixteenth Australasian Workshop on Combinatorial Algorithms, Ballarat, 18-21 September 2005} (Ryan, J., Manyem, P., Sugeng, K. and Miller, M., eds.), University of Ballarat.

\bibitem{Rivest}
R.L. Rivest, A.R. Mayer, D. Kleitman, K. Winklemann, and J. Spencer, ``Coping with errors in binary search procedures'', \textsl{Journal of Computer and System Sciences}, 20, 2, (1980), 396-404.

\bibitem{Minc}
Henryk Minc, \textsl{Nonnegative Matrices}, John Wiley \& Sons, New York, 1988

\bibitem{Lancaster}
Peter Lancaster and Miron Tismenetsky, \textsl{The Theory of Matrices}, second edition, Academic Press, New York, 1985

%\bibitem{Owen}
%Guillermo Owen, \textsl{Game Theory}, Academic Press, New York, 1995

%\bibitem{Shapely}
%L. S. Shapely, ``Stochastic Games'', \textsl{Proceedings of the National Academy of Sciences}, 39, 1095-110

\end{thebibliography}
\end{document}